\documentclass[12pt, showkeys]{article}
\usepackage{amsmath, amssymb, ascmac}
\usepackage{bm}
\usepackage[abbrev]{amsrefs}
\usepackage{enumerate}
\usepackage{amsthm}
\usepackage{comment}
\usepackage[title]{appendix}
\usepackage{color}
\usepackage[dvipdfmx]{graphicx}
\usepackage[subrefformat=parens]{subcaption}

\usepackage{hyperref}

\usepackage[top=30truemm, bottom=30truemm, left=20truemm, right=20truemm]{geometry}

\newtheorem{theorem}{\bf Theorem}
\newtheorem{lem}{\bf Lemma}[section]

\newtheorem{df}{\rm DEFINITION}

\newtheorem{prop}{\bf Proposition}[section]

\title{{\huge Higher order integration by parts formulae \\ for Wiener measures on a path space \\
between two curves}}
\author{{\large Kensuke Ishitani and Soma Nishino}}
\date{}
\begin{document}
\maketitle

\begin{abstract}
We have formulated higher-order integration by parts formulae on the path space restricted between two curves, with respect to pinned/ordinary Wiener measures.
The higher-order integration by parts formulae introduce nontrivial boundary terms, unlike the first-order one. 
Furthermore, in the process of proving these formulae, it becomes necessary to employ the construction methods of Brownian excursion and Brownian house-moving through random walk approximations.
To express the integration by parts formula concisely, we introduced a notation called \emph{Symmetrization.} 
This notation enables the rewriting of the intricate expressions of boundary terms associated with higher-order integration by parts into more concise forms.
Additionally, we provided a probabilistic explanation for the boundary terms by introducing symbols based on the concept of infinitesimal probability. 
These efforts are aimed at fostering an intuitive understanding of the integration by parts formulae.

\end{abstract}


\section{Introduction}
The objective of this study is to provide higher-order integration by parts formulae (IbPF) for the Wiener measure on the path space $K$ constrained between two curves $f^- < f^+$.
The problem of this IbPF is presented in the context of motivating research on the reflected stochastic partial differential equations (SPDEs) by Zambotti.
The IbPF for the first-order differential case has been proved by Zambotti~\cite{zambotti} in the scenario where $f^- = \text{constant}$ and $f^+ = \infty$.
For the general case of $f^-, f^+ \in W^{1,2+}([0,1])$, Funaki and Ishitani~\cite{funaki_ishitani} have established the integration by parts formulae.
According to Funaki and Ishitani, the following integration by parts formula holds:
\begin{equation}\nonumber
\begin{aligned}
    &\quad \int_{\Omega} 1_K(X^{a,b}) \nabla^{(1)} \varphi (h) (X^{a,b}) \ \mathrm{d} P \\
    &\quad = \int_{\Omega} 1_K(X^{a,b}) \varphi(x) \int_{0}^{1} h'(t) \ \mathrm{d} X^{a,b}(t) \  \mathrm{d} P \\
    &\qquad + \int_{0}^{1} h(t) \nu_{(1)}^+(t) \int_{\Omega} \varphi \left( X_{[0,t]}^{a,f^+(t),(f^-,f^+)} \oplus X_{[t,1]}^{f^+(t),b,(f^-,f^+)} \right) \ \mathrm{d} P \ \mathrm{d} t \\
    &\qquad - \int_{0}^{1} h(t) \nu_{(1)}^-(t) \int_{\Omega} \varphi \left( X_{[0,t]}^{a,f^-(t),(f^-,f^+)} \oplus X_{[t,1]}^{f^-(t),b,(f^-,f^+)} \right) \ \mathrm{d} P \ \mathrm{d} t .
\end{aligned}
\end{equation}

When extending this formula to the case of $d$-th order differentials, simplifying the description is crucial due to the complexity of calculations.
We have successfully achieved this by introducing appropriate symbols. Moreover, by incorporating the concept of infinitesimal probability, we give a probabilistic meaning to $\nu^{\pm}$ appearing in the boundary terms of the integration by parts formulae.
Additionally, employing a notation similar to the Symmetrization notation in the theory of multiple Wiener integrals~\cite{nualart}, we have succinctly described the formulae for the $d$-th order differentials.

Furthermore, in the case of $d$ ($\geq 2$)-th order differentials, new construction methods involving random walk approximations of Brownian excursion and Brownian house-moving~\cite{yanashima} become necessary during the proof of the formulae.
Considering these formulae for Wiener measures on the time interval $[0,1]$, we have provided integration by parts formulae for two cases: one where the value is pinned at $t=1$ and the other where it is free.

\section{Notation}
For $0 \leq s < t \leq 1$, let $C([s,t], \mathbb{R})$ be the class of $\mathbb{R}$-valued continuous functions defined on $[s,t]$.
For $w_1, w_2 \in C([s,t], \mathbb{R})$, we define the metric $d_{\infty}$ by
\begin{equation}\nonumber
d_{\infty}(w_1,w_2) = \sup_{u \in [s,t]} |w_1(u) - w_2(u)|.
\end{equation}
$\mathcal{B}(C([s,t], \mathbb{R}))$ denotes the Borel $\sigma$-algebra with respect to the topology generated by the metric $d_{\infty}$.

Suppose that $Y ~ \colon ~ (\Omega, \mathcal{F}, P) \to (C([0,1], \mathbb{R}), \mathcal{B}(C([0,1], \mathbb{R})))$ is a random variable and that $\Lambda \in \mathcal{B}(C([0,1], \mathbb{R}))$ satisfies $P(Y \in \Lambda) > 0$.
Then, we define the probability measure $P_{Y^{-1}(\Lambda)}$ on $(Y^{-1}(\Lambda), Y^{-1}(\Lambda) \cap \mathcal{F})$ as
\begin{equation}\nonumber
    P_{Y^{-1}(\Lambda)}(A) := \frac{P(A)}{P(Y \in \Lambda)}, ~ A \in Y^{-1}(\Lambda) \cap \mathcal{F} := \{ Y^{-1}(\Lambda) \cap F \mid F \in \mathcal{F} \} .
\end{equation}
Let $Y\vert_{\Lambda}$ denote the restriction $Y$ to $(Y^{-1}(\Lambda), Y^{-1}(\Lambda) \cap \mathcal{F}, P_{Y^{-1}(\Lambda)})$. Then,
\begin{equation}\nonumber
Y\vert_{\Lambda} ~ \colon ~ (Y^{-1}(\Lambda), Y^{-1}(\Lambda) \cap \mathcal{F}, P_{Y^{-1}(\Lambda)}) \to (\Lambda, \mathcal{B}(\Lambda)),
\end{equation}
is a random variable.
In this study, we write $P_{Y^{-1}(\Lambda)}(Y\vert_{\Lambda} \in \Gamma)$ as $P(Y\vert_{\Lambda} \in \Gamma)$ and $E^{P_{Y^{-1}(\Lambda)}}[f(Y\vert_{\Lambda})]$ as $E[f(Y\vert_{\Lambda})]$ respectively.

For $t > 0$ and $x, y \in \mathbb{R}$, we define the heat kernel $p(t;x,y)$ as
\begin{equation}\nonumber
p(t;x,y) := \frac{1}{\sqrt{2\pi t}} \exp \left( -\frac{(x-y)^2}{2t} \right).
\end{equation}

For $0 \leq t_1 < t_2 \leq 1$ and $f, g \in C([0,1], \mathbb{R})$, we define
\begin{equation}\nonumber
    K_{[t_1,t_2]}(f,g) := \{ w = \{ w(t) \}_{t \in [t_1,t_2]} \in C([t_1,t_2], \mathbb{R}) \mid f(t) \leq w(t) \leq g(t), ~ t_1 \leq t \leq t_2 \} ,
\end{equation}
\begin{equation}\nonumber
K_{[t_1,t_2]}^+(f) := \bigcup_{n=1}^{\infty} K_{[t_1,t_2]}(f,n), ~ K_{[t_1,t_2]}^-(f) := \bigcup_{n=1}^{\infty} K_{[t_1,t_2]}(-n,f).
\end{equation}
Moreover, we write simply
\begin{equation}\nonumber
K(f,g) := K_{[0,1]}(f,g), ~ K^+(f) := K_{[0,1]}^+(f), ~ K^-(f) := K_{[0,1]}^-(f) .
\end{equation}

Let $W_{[t_1,t_2]} = \{ W(t) \}_{t \in [t_1,t_2]}$ and $B_{[t_1,t_2]}^{a \to b} = \{ B^{a \to b}(t) \}_{t \in [t_1,t_2]}$ denote the standard one-dimensional Brownian motion on $[t_1,t_2]$, and the one-dimensional Brownian bridge from $a$ to $b$ on $[t_1, t_2]$, respectively.
Let $f^-,\ f^+$ be $\mathbb{R}$-valued $C^2$-functions on $[0,1]$, satisfying the condition:
\begin{equation}\label{twocurves}
    \min_{0 \leq t \leq 1} (f^+(t) - f^-(t)) > 0 .
\end{equation}
Here, for $0 \leq t_1 < t_2 \leq 1$ and $f^-(t_1) \leq a \leq f^+(t_1)$, $f^-(t_2) \leq b \leq f^+(t_2)$, we define a continuous stochastic process $X_{[t_1,t_2]}^{a,b,(f^-,f^+)}$ on $[t_1,t_2]$ as follows:
\begin{description}
    \item[{\rm (i)}] in the case $f^-(t_1) < a < f^+(t_1)$ and $f^-(t_2) < b < f^+(t_2)$, the conditioned process $B_{[t_1,t_2]}^{a \to b}\vert_{K_{[t_1,t_2]}(f^-,f^+)}$;
    \item[{\rm (ii)}] in the case $a = f^{\pm}(t_1)$ and $f^-(t_2) < b < f^+(t_2)$, the weak limit of $B_{[t_1,t_2]}^{a \to b}\vert_{K_{[t_1,t_2]}(f^- -\varepsilon,f^+ + \varepsilon)}$ as $\varepsilon \downarrow 0$;
    \item[{\rm (iii)}] in the case $f^-(t_1) < a < f^+(t_1)$ and $b = f^{\pm}(t_2)$, the weak limit of $B_{[t_1,t_2]}^{a \to b}\vert_{K_{[t_1,t_2]}(f^- -\varepsilon,f^+ + \varepsilon)}$ as $\varepsilon \downarrow 0$;
    \item[{\rm (iv)}] in the case $a = f^{\pm}(t_1)$ and $b = f^{\pm}(t_2)$, the weak limit of $B_{[t_1,t_2]}^{a \to b}\vert_{K_{[t_1,t_2]}(f^- -\varepsilon,f^+ + \varepsilon)}$ as $\varepsilon \downarrow 0$;
    \item[{\rm (v)}]in the case $a = f^{\pm}(t_1)$ and $b = f^{\mp}(t_2)$, the weak limit of $B_{[t_1,t_2]}^{a \to b}\vert_{K_{[t_1,t_2]}(f^- -\varepsilon,f^+ + \varepsilon)}$ as $\varepsilon \downarrow 0$.
\end{description}
In this definition, in cases (\textrm{ii}) and (\textrm{iii}), the weak limits have the same probability law as a three-dimensional Bessel bridge between the two curves $f^+$ and $f^-$.
In case (\textrm{iv}), the weak limit is a Brownian excursion, and in case (\textrm{v}), the weak limit is a continuous Markov process called a Brownian house-moving~\cite{ishitani_hatakenaka_suzuki}.
Furthermore, for $f^-(t_1) \leq a \leq f^+(t_1)$, we define a continuous stochastic process $X_{[t_1,t_2]}^{a,(f^-,f^+)}$ on $[t_1,t_2]$ as follows:
\begin{description}
    \item[{\rm (vi)}] in the case $f^-(t_1) < a < f^+(t_1)$, the conditioned process $(a + W_{[t_1,t_2]})\vert_{K_{[t_1,t_2]}(f^-,f^+)}$;
    \item[{\rm (vii)}] in the case $a = f^{\pm}(t_1)$, the weak limit of $(a + W_{[t_1,t_2]})\vert_{K_{[t_1,t_2]}(f^- - \varepsilon, f^+ + \varepsilon)}$ as $\varepsilon \downarrow 0$.
\end{description}
In the case (\textrm{vii}), the weak limit is a Brownian meander.
Additionally, let $X_{[t_1,t_2]}^{a,b} \equiv B_{[t_1,t_2]}^{a \to b}$ and $X_{[t_1,t_2]}^{a} \equiv a + W_{[t_1,t_2]}$ denote the Brownian bridge from $a$ to $b$ and the Brownian motion starting from $a$ for consistency in the terminology.
For $\mathbb{R}$-valued continuous function $X$ and $C^2$-function $g$ on $[t_1,t_2]$, we define the Cameron--Martin density as
\begin{equation}\nonumber
    \begin{aligned}
        &Z_{[t_1,t_2]}^{g}(X) \\
        &:= \exp \left( g'(t_2) X(t_2) - g'(t_1) X(t_1) - \int_{t_1}^{t_2} X(u) g''(u) \ \mathrm{d} u - \frac{1}{2} \int_{t_1}^{t_2} (g'(u))^2 \ \mathrm{d} u \right).
    \end{aligned}
\end{equation}

Let $t_0 < t_1 < t_2$.
Suppose that $X_{[t_0,t_1]}^{(1)} \in C([t_0,t_1], \mathbb{R})$ and $X_{[t_1,t_2]}^{(2)} \in C([t_1,t_2], \mathbb{R})$ satisfy $X_{[t_0,t_1]}^{(1)}(t_1) = X_{[t_1,t_2]}^{(2)}(t_1)$.
In this case, we set
\begin{equation}\nonumber
(X_{[t_0,t_1]}^{(1)} \oplus X_{[t_1,t_2]}^{(2)} )(t) :=
\begin{cases}
X_{[t_0,t_1]}^{(1)}(t), & t_0 \leq t \leq t_1 , \\
X_{[t_1,t_2]}^{(2)}(t), & t_1 \leq t \leq t_2 .
\end{cases}
\end{equation}
Let $t_0 < t_1 < \dots < t_j$. For $X_{[t_{i-1},t_i]}^{(i)} \in C([t_{i-1},t_i], \mathbb{R})$ $(i = 1, 2, \dots, j)$ that satisfy $X_{[t_{i-1},t_i]}^{(i)}(t_i) = X_{[t_{i},t_{i+1}]}^{(i+1)}(t_i)$ $(i = 1, 2, \dots, j-1)$, we define
\begin{equation}\nonumber
    \bigoplus_{i=1}^j X_{[t_{i-1},t_i]}^{(i)} (t) := X_{[t_{k-1},t_k]}^{(k)}(t) ,~ t_{k-1} \leq t \leq t_k .
\end{equation}

As the set of time parameters, we use
\begin{equation}\nonumber
    \mathcal{T}^{[m]} := \{ t = (t_j)_{j=1}^m \in [0,1]^m \mid t_1 < t_2 < \cdots < t_m \} .
\end{equation}
In addition, we often write the closed interval $[0,1]$ as $\mathcal{T}^{[1]}$ for consistency in the terminology.

Let $W^{1,p}([0,1])$ be the Sobolev space on $[0,1]$.
We define $W^{1,2+}([0,1]) := \bigcup_{p>2} W^{1,p}([0,1])$.
Let $H$ be the Hilbert space $L^2([0,1])$ equipped with the standard inner product $\langle \bullet, \bullet \rangle_H$, and let $H_0^1 \equiv H_0^1((0,1))$ denote the closure of $C_c^{\infty}((0,1))$ in $H^1([0,1]) \equiv W^{1,2}([0,1])$.
We define $C_b^d(H)$ as the set of $d$-th $H_0^1$ differentiable function on $H$, which is bounded and has bounded continuous $j$-th $H_0^1$ derivatives~($j = 1, \dots, d$).

For $\varphi \in C_b^d(H)$, let $\nabla^{(d)} \varphi$ be the $d$-th $H_0^1$ derivative of $\varphi$.
The definition of $H_0^1$ derivatives is given below.
\begin{df}
Let $\varphi ~ \colon ~ H \to \mathbb{R}$ be $d$-th $H_0^1$ differentiable at the point $x \in H$.
This means that there exists a continuous $d$-linear map $\Phi \colon H \times \dots \times H \to \mathbb{R}$ such that, for any $h_1, \dots, h_d \in H_0^1$, the following equality holds:
\begin{equation}\label{H01}
    \frac{\partial^d \varphi}{\partial t_1 \cdots \partial t_d}(x + t_1 h_1 + \dots + t_d h_d) \vert_{t_1 = \dots = t_d = 0} = \Phi (h_1, \dots, h_d) .
\end{equation}
In this case, we define $\Phi$ as the $d$-th order $H_0^1$ derivative of $\varphi$ at the point $x$, denoted as $\nabla^{(d)} \varphi$.
Throughout this paper, the right-hand side of equation (\ref{H01}) is often written as $\nabla^{(d)} \varphi(h_1, \dots, h_d)(x)$.
\end{df}

In this study, we assume that the following condition holds for $\varphi \in C_b^d(H)$:
\begin{equation}\label{phi_cond}
    \exists \ell \in \mathbb{N}, \exists \Phi \in C_b^d(\mathbb{R}^{\ell}), \exists \lambda_1, \dots, \lambda_{\ell} \in L^{d+}([0,1]) ~\text{s.t.}~ \varphi(x) = \Phi(\langle x, \lambda_1 \rangle_H , \dots, \langle x, \lambda_{\ell} \rangle_H ) .
\end{equation}

\section{Main results}
\subsection{The case with two pinned edges}
\begin{theorem}\label{main1}
Assume that $f^-,\ f^+ \in W^{1,2+}([0,1])$, $f^-(0) < a < f^+(0),\ f^-(1) < b < f^+(1)$, and that condition (\ref{twocurves}) holds for the two curves $f^-$ and $f^+$.
Then, for every $h_1, h_2, \dots, h_d \in H_0^1((0,1))$ and $\varphi \in C_b^d(H)$ satisfying the condition (\ref{phi_cond}), we have the integration by parts formula
\begin{equation}\nonumber
    \begin{aligned}
        &\int_{\Omega} 1_K(X^{a,b}) \nabla^{(d)} \varphi(h_1,h_2, \dots , h_d)(X^{a,b}) \ \mathrm{d} P \\
        &= \int_{\Omega} 1_K(X^{a,b}) \varphi(X^{a,b}) \left( \prod_{i=1}^{d} \int_{\mathcal{T}^{[1]}} h_i'(t_i) \ \mathrm{d} X^{a,b}(t_i) \right) \ \mathrm{d} P \\
        &\quad + \text{BD}^{(1)} + \text{BD}^{(2)} + \dots + \text{BD}^{(d)}.
    \end{aligned}
\end{equation}
Here, the boundary term is
\begin{equation}\nonumber
    \begin{aligned}
        \text{BD}^{(j)}&:= \sum_{\varepsilon_1 , \dots , \varepsilon_j \in \{ 1, -1 \} } \varepsilon_1 \cdots \varepsilon_j \int_{\mathcal{T}^{[j]}} \frac{1}{(d-j)!} \nu_{(j)}^{\varepsilon_1 , \dots , \varepsilon_j}(t_1, \dots , t_j) \\
        &\quad \times \int_{\Omega} \varphi(Y_{t_1, \dots , t_j}^{\varepsilon_1, \dots, \varepsilon_j, (f^-,f^+)}) \sum_{\sigma \in \mathfrak{S}_d} h_{\sigma(1)}(t_1) \cdots h_{\sigma(j)}(t_j) \\
        &\quad \times \left( \prod_{i=j+1}^{d} \int_{\mathcal{T}^{[1]}} h_{\sigma(i)}'(t_{i}) ~ \ \mathrm{d} Y_{t_1, \dots , t_j}^{\varepsilon_1, \dots, \varepsilon_j, (f^-,f^+)} (t_{i}) \right) ~ \ \mathrm{d} P \ \mathrm{d} t_1 \cdots \ \mathrm{d} t_j.
    \end{aligned}
\end{equation}
The following is an explanation of the symbols used in the statement.
Now, we have denoted $f^{\pm 1} := f^{\pm}$.
\begin{equation}\nonumber
    \begin{aligned}
        &Y_{t_1, \dots , t_j}^{\varepsilon_1, \dots, \varepsilon_j, (f^-,f^+)} \\
        &:= X_{[0,t_1]}^{a,f^{\varepsilon_1}(t_1),(f^-,f^+)} \oplus \left( \bigoplus_{i=2}^{j} X_{[t_{i-1},t_i]}^{f^{\varepsilon_{i-1}}(t_{i-1}),f^{\varepsilon_i}(t_i),(f^-,f^+)} \right) \oplus X_{[t_j,1]}^{f^{\varepsilon_j}(t_j),b,(f^-,f^+)} ,
    \end{aligned}
\end{equation}

\begin{equation}\nonumber
    \begin{aligned}
        &\, \, \nu_{(j)}^{\varepsilon_1, \dots, \varepsilon_j}(t_1, \dots, t_j) \\
        &:= p^{(j)}_{a,b}(t_1, \dots, t_j ; f^{\varepsilon_1}(t_1), \dots, f^{\varepsilon_j}(t_j)) \\
        &\quad \times \frac{\Delta P_{[0,t_1]}^{a,f^{\varepsilon_1}(t_1)} (K_{[0,t_1]})}{\sqrt{t_1}} \left( \prod_{i=2}^{j} \frac{\Delta^2 P_{[t_{i-1},t_i]}^{f^{\varepsilon_{i-1}}(t_{i-1}),f^{\varepsilon_i}(t_i)} (K_{[t_{i-1},t_i]})}{(\sqrt{t_{i} - t_{i-1}})^2} \right) \frac{\Delta P_{[t_j,1]}^{f^{\varepsilon_j}(t_j), b} (K_{[t_j,1]})}{\sqrt{1-t_j}}
    \end{aligned}
\end{equation}
For $0 = t_0 < t_1 < t_2 < \cdots < t_j < 1$ and $c_0 = a, c_1, c_2 , \dots , c_j \in \mathbb{R}$,
\begin{equation}\nonumber
    \begin{aligned}
        &\, p^{(j)}_{a,b}(t_1, \dots, t_j ; c_1 , \dots, c_j) \\
        &:= \prod_{i=1}^{j} p(t_i - t_{i-1}; c_{i-1}, c_i) \frac{p(1-t_j; c_j, b)}{p(1;a,b)} \\
        &= \frac{P(X^{a,b}(t_1) \in \ \mathrm{d} c_1)}{\ \mathrm{d} c_1} \prod_{i=2}^{j} \frac{P(X^{a,b}(t_i) \in \ \mathrm{d} c_i \mid X^{a,b}(t_{i-1}) = c_{i-1})}{\ \mathrm{d} c_i}.
    \end{aligned}
\end{equation}
For these symbols, especially $p^{(j)}_{a,b}$ is the finite-dimensional density of the Brownian bridge, and $\Delta P$ and $\Delta^2 P$ can be interpreted as first and second-order infinitesimal probabilities, respectively.
That is, for the $\nu_{(d)}^{\varepsilon_1, \dots, \varepsilon_j}$, which was unclear in meaning in~\cite{funaki_ishitani}, a probabilistic interpretation becomes possible.
Furthermore, the expression $\Delta P / \sqrt{t}$ can be interpreted as the scaling of the first-order infinitesimal probability.
The same holds for the second-order case.
The specific expressions for infinitesimal probabilities will be presented in a later lemma \ref{infinitesimal_lemma}.
\end{theorem}

\subsection{The cases with edges one pinned and the other free}
\begin{theorem}\label{main2}
Assume that $f^-,\ f^+ \in W^{1,2+}([0,1])$, $f^-(0) < a < f^+(0)$, and that condition (\ref{twocurves}) holds for the two curves $f^-$ and $f^+$.
Then, for every $h_1, h_2, \dots, h_d \in H_0^1((0,1))$ and $\varphi \in C_b^d(H)$ satisfying the condition (\ref{phi_cond}), we have the integration by parts formula
\begin{equation}\nonumber
    \begin{aligned}
        &\int_{\Omega} 1_K(X^{a}) \nabla^{(d)} \varphi(h_1,h_2, \dots , h_d)(X^{a}) \ \mathrm{d} P \\
        &= \int_{\Omega} 1_K(X^{a}) \varphi(X^{a}) \left( \prod_{i=1}^{d} \int_{\mathcal{T}^{[1]}} h_d'(t_d) \ \mathrm{d} X^{a}(t_d) \right) \ \mathrm{d} P \\
        &\quad + \text{BD}^{(1)} + \text{BD}^{(2)} + \dots + \text{BD}^{(d)}.
    \end{aligned}
\end{equation}
The symbols used in the assertions are explained below.
Here, we have denoted $f^{\pm 1} := f^{\pm}$.
\begin{equation}\nonumber
    \begin{aligned}
        &\, Y_{t_1, \dots , t_j}^{\varepsilon_1, \dots, \varepsilon_j, (f^-,f^+)} \\
        &:= X_{[0,t_1]}^{a,f^{\varepsilon_1}(t_1),(f^-,f^+)} \oplus \left( \bigoplus_{i=2}^{j} X_{[t_{i-1},t_i]}^{f^{\varepsilon_{i-1}}(t_{i-1}),f^{\varepsilon_i}(t_i),(f^-,f^+)} \right) \oplus X_{[t_j,1]}^{f^{\varepsilon_j}(t_j),(f^-,f^+)} ,
    \end{aligned}
\end{equation}

\begin{equation}\nonumber
    \begin{aligned}
        &\, \, \nu_{(j)}^{\varepsilon_1, \dots, \varepsilon_j}(t_1, \dots, t_j) \\
        &:= p^{(j)}_{a}(t_1, \dots, t_j; f^{\varepsilon_1}(t_1), \dots, f^{\varepsilon_j}(t_j)) \\
        &\quad \times \frac{\Delta P_{[0,t_1]}^{a,f^{\varepsilon_1}(t_1)} (K_{[0,t_1]})}{\sqrt{t_1}} \left( \prod_{i=2}^{j} \frac{\Delta^2 P_{[t_{i-1},t_i]}^{f^{\varepsilon_{i-1}}(t_{i-1}),f^{\varepsilon_i}(t_i)} (K_{[t_{i-1},t_i]})}{(\sqrt{t_{i} - t_{i-1}})^2} \right) \frac{\Delta P_{[t_j,1]}^{f^{\varepsilon_j}(t_j)} (K_{[t_j,1]})}{\sqrt{1-t_j}} .
    \end{aligned}
\end{equation}
For $0 = t_0 < t_1 < t_2 < \cdots < t_j < 1$ and $c_0 = a, c_1, c_2 , \dots , c_j \in \mathbb{R}$,
\begin{equation}\nonumber
    \begin{aligned}
        &\, p^{(j)}_{a}(t_1, \dots, t_j ; c_1 , \dots, c_j) \\
        &:= \prod_{i=1}^{j} p(t_i - t_{i-1}; c_{i-1}, c_i) \\
        &= \frac{P(X^{a}(t_1) \in \ \mathrm{d} c_1)}{\ \mathrm{d} c_1} \prod_{i=2}^{j} \frac{P(X^{a}(t_i) \in \ \mathrm{d} c_i \mid X^{a}(t_{i-1}) = c_{i-1})}{\ \mathrm{d} c_i}.
    \end{aligned}
\end{equation}
With respect to these symbols, $p^{(j)}_{a}$ is the finite-dimensional density of Brownian motion unlike $p^{(j)}_{a,b}$.
\end{theorem}

\subsection{Related results}
Zambotti~\cite{zambotti} established the first-order integration by parts formula on the restricted path space $K(f^-,f^+)$ in the case of $f^- = \text{constant}$ and $f^+ = \infty$, applying it to the solutions of Nualart and Pardoux type SPDE. Funaki--Ishitani~\cite{funaki_ishitani} extended this to the general case of $f^-, f^+ \in W^{1,2+}([0,1])$.

Our proof technique aligns with that used by Funaki--Ishitani~\cite{funaki_ishitani}, relying on the classical polygonal approximation for Brownian motion.
It reduces the path-space integration to integration in finite-dimensional spaces.
Additionally, alternative approaches utilize Denisov's formula and Biane's formula for Brownian motion~\cites{bonaccorsi,zambotti}, or use hitting times~\cite{hariya}.
In comparison, our method is straightforward and advantageous as it can be applied to a wide class of domains $K$.

The integration by parts formulae in infinite-dimensional spaces have been discussed within the framework of abstract Wiener spaces by Goodman~\cite{goodman} and Shigekawa~\cite{shigekawa}.
Relevant to our research, Otobe~\cite{Otobe} investigated solutions restricted between two curves for the SPDEs of Nualart and Pardoux type.
In addition, Hariya~\cite{hariya} studied the first-order integration by parts formula on Wiener space for paths taking values in $\mathbb{R}^d$.
Debussche--Zambotti~\cite{debussche_zambotti} and Debussche--Goudenege~\cite{debussche_goudenege} explored solutions to the Stochastic Cahn--Hilliard equation.
Altman--Zambotti~\cite{altman_zambotti} analyzed the stochastic partial differential equations associated with Bessel bridges.
Moreover, Zambotti~\cite{zambotti2} provided a comprehensive overview of a series of topics related to the integration by parts formulae.

\section{Polygonal approximation}
For a finite partition $\wp = \{0 = t_0 < t_1 < \dots < t_m = 1\}$ of the interval $[0,1]$, let $\pi_{\wp} ~ \colon ~ C([0,1],\mathbb{R}) \to C([0,1],\mathbb{R})$ be the polygonalization of $x \in C([0,1],\mathbb{R})$ with respect to $\wp$.
In other words, $\pi_{\wp} = \pi_{\wp,2} \circ \pi_{\wp,1}$, where $\pi_{\wp,1} ~ \colon ~ C([0,1],\mathbb{R}) \to \mathbb{R}^{m+1}$ is given by $\pi_{\wp,1}x \equiv ((\pi_{\wp,1}x)_k)_{k=0}^{m} = (x(t_k))_{k=0}^{m}$, and $\pi_{\wp,2} ~ \colon ~ \mathbb{R}^{m+1} \to C([0,1],\mathbb{R})$ is defined for $\underline{x} = (x_k)_{k=0}^{m} \in \mathbb{R}^{m+1}$ as
\begin{equation}\nonumber
    (\pi_{\wp,2} \underline{x})(t) := \frac{(t-t_{k-1})x_k + (t_k - t)x_{k-1}}{t_k - t_{k-1}}, \quad t \in [t_{k-1},t_k]
\end{equation}
for $1 \leq k \leq m$.
The mesh size of this finite partition is defined as $|\wp| = \max_{1\leq k \leq m}(t_k - t_{k-1})$.

Now, suppose we are given a sequence $\{ \wp_n \}_{n=1}^{\infty}$ of finite partitions of the interval $[0,1]$ such that $\lim_{n \to \infty} |\wp_n| = 0$.
Denote $\pi_n = \pi_{\wp_n}$ and $\pi_{n,i} = \pi_{\wp_n,i}$ for $i=1,2$.
In this case, the integration in the path space can be approximated by integration in a finite-dimensional space.
Specifically, the following Proposition holds~\cite{funaki_ishitani}.
\begin{prop}\label{polygonal}
    For every probability measure $P$ on $C([0,1],\mathbb{R})$ and bounded continuous function $\Phi$ on $C([0,1],\mathbb{R})$, we have
    \begin{equation}\nonumber
        \begin{aligned}
            E^P[\Phi;K] &= \lim_{n \to \infty} E^P[\Phi \circ \pi_n ; \pi_n^{-1}(K(\pi_n f^-, \pi_n f^+))] \\
            &= \lim_{n \to \infty} E^{P \circ \pi_{n,1}^{-1}}[\Phi \circ \pi_{n,2}; \underline{K}_{\wp_n}].
        \end{aligned}
    \end{equation}
\end{prop}

\section{Proof of Theorem \ref{main1}}
Let $f^{\pm} \in W^{1,p}([0,1])$, where $p > 2$.
We can assume that $f_n^{\pm}(0) = f^{\pm}(0)$ and $f_n^{\pm}(1) = f^{\pm}(1)$ hold for a sequence $\{ f_n^{\pm} \}_{n=1}^{\infty} \subset C^2([0,1])$ that approximates $f^{\pm}$ in the Sobolev norm.
Therefore, it is sufficient to establish the case when $f^{\pm} \in C^2([0,1])$.
In this case, $\nu_{(j)}^{\varepsilon_1, \dots, \varepsilon_j}$ converges, as shown in appendix \ref{smooth_curves}.

Let $h \in C([0,1], \mathbb{R})$ be a function with compact support within $(0,1)$, denoted by $C_c((0,1))$, and we define $C_c^2((0,1)) := C^2([0,1]) \cap C_c((0,1))$.
The integration by parts formula in Theorem~\ref{main1} is continuous under the $H^1$ norm. Hence, it suffices to establish the integration by parts formula for $h_1, h_2, \dots , h_d \in C_c^2((0,1))$.

In the proof of this proposition, it is convenient to consider the continuous function space $C([0,1])$ instead of the general sample space $\Omega$.
As index space symbols, we use
\begin{equation}\nonumber
    \mathcal{K}_{n-1}^{[m]} := \{ k = (k_j)_{j=1}^m \in \{ 1, 2, \cdots, n-1 \}^m \mid k_1 < k_2 < \cdots < k_m \} .
\end{equation}
For consistency in the terminology, we write the index space $\{ 1, 2, \dots, n-1 \}$ as $\mathcal{K}_{n-1}^{[1]}$.
To simplify later discussions, we partition the interval $[0,1]$ into intervals of width $1/n$: $\wp_n = \{ t_k = k/n \}_{k=0}^{n}$.
By applying Proposition~\ref{polygonal} to the function $\nabla^{(d)}\varphi(h_1, h_2, \dots, h_d)$ and the pinned Wiener measure $P^{a,b} \equiv P(X^{a,b})^{-1}$, we obtain the following equality:
\begin{equation}\nonumber
    \begin{aligned}
        &\int_{\Omega} 1_K(X^{a,b}) \nabla^{(d)} \varphi(h_1,h_2, \dots , h_d)(X^{a,b}) \ \mathrm{d} P \\
        &= \lim_{n \to \infty} \int_{\underline{K}_n} \sum_{k_1, \dots , k_d \in \mathcal{K}_{n-1}^{[1]}} \frac{\partial^d \tilde{\varphi}_n}{\partial x_{k_1} \cdots \partial x_{k_d}} (\underline{x}) (\pi_{n,1} h_1)_{k_1} \cdots (\pi_{n,1} h_d)_{k_d} P^{a,b} \circ \pi_{n,1}^{-1}(\ \mathrm{d} \underline{x})
    \end{aligned}
\end{equation}
Here, $\underline{K}_n \equiv \underline{K}_{\wp_n}(f^-,f^+) = \prod_{j=0}^{n} [f^-(t_j), f^+(t_j)] \subset \mathbb{R}^{n+1}$, $\tilde{\varphi}_n(\underline{x}) = \varphi(\pi_{n,2} \underline{x})$, and $\underline{x} = (x_j)_{j=0}^{n} \in \mathbb{R}^{n+1}$ with $\ \mathrm{d} \underline{x} = \prod_{j=0}^{n} \ \mathrm{d} x_j$.
Note that $(\pi_{n,1}h_j)_{k_j} \equiv h_j(t_{k_j}) = h_j(k_{j}/n)$.
Indeed,
\begin{equation}\nonumber
    \begin{aligned}
        &\quad |\nabla^{(d)}\varphi(h_1, h_2, \dots, h_d) - \nabla^{(d)}\varphi(\pi_n h_1, \pi_n h_2, \dots, \pi_n h_d)| \\
        &\lesssim ||\nabla^{(d)} \varphi||(||h_1 - \pi_n h_1||_H + ||h_2 - \pi_n h_2||_H + \dots + ||h_d - \pi_n h_d||_H) \\
        &\to 0 ~ (n \to \infty),
    \end{aligned}
\end{equation}
so this finite-dimensional approximation holds in the same way for Funaki--Ishitani~\cite{funaki_ishitani}.
Here, $||\nabla^{(d)} \varphi||$ is the operator norm of the $d$-th order $H_0^1$ derivative $\nabla^{(d)} \varphi$.
Moreover, the supports of $h_1, h_2, \dots, h_d$ are compact within $(0,1)$, so terms corresponding to $k=0$ and $k=n$ do not appear in the sum.

For $0 \leq k_1 < k_2 \leq n$ and $\underline{x} = (x_j)_{j=k_1}^{k_2} \in \mathbb{R}^{k_2 - k_1 + 1}$, we define
\begin{equation}\nonumber
    q_n^{k_1,k_2}(\underline{x}) := \exp \left( -\frac{n}{2} \sum_{j=k_1 +1}^{k_2} (x_j - x_{j-1})^2 \right),
\end{equation}
and for $1 \leq k \leq n$, we set
\begin{equation}\nonumber
    \Xi_{n,a,b}^{k} := \frac{1}{\sqrt{k}} \left( \frac{2\pi}{n} \right)^{\frac{k-1}{2}} \exp \left( -\frac{n (a-b)^2}{2k} \right).
\end{equation}
Then, we have
\begin{equation}\nonumber
    \begin{aligned}
        P^{a,b} \circ \pi_{n,1}^{-1}(\ \mathrm{d} \underline{x}) &= \frac{1}{p(1;a,b)} \prod_{j=1}^{n} p(1/n;x_{j-1},x_j) \ \mathrm{d} \underline{x}^{a,b} \\
        &= (\Xi_{n,a,b}^n)^{-1} q_n^{0,n}(\underline{x}) \ \mathrm{d} \underline{x}^{a,b}.
    \end{aligned}
\end{equation}
Here, in the right-hand side, for $\underline{x} = (x_j)_{j=0}^n \in \mathbb{R}^{n+1}$, we define $\ \mathrm{d} \underline{x}^{a,b} := \delta_a(\ \mathrm{d} x_0) \prod_{j=1}^{n-1} \ \mathrm{d} x_j \delta_b(\ \mathrm{d} x_n)$.
Therefore, the finite-dimensional approximation formula can be written as
\begin{equation}\nonumber
    \begin{aligned}
        &\int_{\underline{K}_n} \sum_{k_1, \dots , k_d \in \mathcal{K}_{n-1}^{[1]}} \frac{\partial^d \tilde{\varphi}_n}{\partial x_{k_1} \cdots \partial x_{k_d}} (\underline{x}) (\pi_{n,1} h_1)_{k_1} \cdots (\pi_{n,1} h_d)_{k_d} P^{a,b} \circ \pi_{n,1}^{-1}(\ \mathrm{d} \underline{x}) \\
        &= \sum_{k_1, \dots , k_d \in \mathcal{K}_{n-1}^{[1]}} h_1 \left( \frac{k_1}{n} \right) \cdots h_d \left( \frac{k_d}{n} \right) \int_{\underline{K}_n} \frac{\partial^d \tilde{\varphi}_n}{\partial x_{k_1} \cdots \partial x_{k_d}} (\underline{x}) (\Xi_{n,a,b}^n)^{-1} q_n^{0,n}(\underline{x}) \ \mathrm{d} \underline{x}^{a,b}.
    \end{aligned}
\end{equation}
Applying the integration by parts formula on the finite-dimensional space $\mathbb{R}^{n-1}$ to this expression, we obtain
\begin{equation}\nonumber
    \begin{aligned}
        &\sum_{k_1, \dots , k_d \in \mathcal{K}_{n-1}^{[1]}} h_1 \left( \frac{k_1}{n} \right) \cdots h_d \left( \frac{k_d}{n} \right) \int_{\underline{K}_n} \frac{\partial^d \tilde{\varphi}_n}{\partial x_{k_1} \cdots \partial x_{k_d}} (\underline{x}) (\Xi_{n,a,b}^n)^{-1} q_n^{0,n}(\underline{x}) \ \mathrm{d} \underline{x}^{a,b} \\
        &= (-1)^d \sum_{k_1, \dots , k_d \in \mathcal{K}_{n-1}^{[1]}} h_1 \left( \frac{k_1}{n} \right) \cdots h_d \left( \frac{k_d}{n} \right) \int_{\underline{K}_n} \tilde{\varphi}_n(\underline{x}) (\Xi_{n,a,b}^n)^{-1} \frac{\partial^d q_n^{0,n}}{\partial x_{k_d} \cdots \partial x_{k_1}} (\underline{x}) \ \mathrm{d} \underline{x}^{a,b} \\
        &\quad + \text{BD}_n^{(1)} + \text{BD}_n^{(2)} + \dots + \text{BD}_n^{(d)}.
    \end{aligned}
\end{equation}
Here, the boundary terms are given by
\begin{equation}\nonumber
    \begin{aligned}
        &\text{BD}_n^{(j)} \\
        &:= \sum_{\varepsilon_1 , \dots , \varepsilon_j \in \{ 1, -1 \}} (-1)^{d-j} \varepsilon_1 \cdots \varepsilon_j \sum_{\substack{(k_1, \dots , k_j) \in \mathcal{K}_{n-1}^{[j]} \\ k_{j+1}, \dots , k_d \in \mathcal{K}_{n-1}^{[1]}}} \frac{1}{(d-j)!} \\
        &\quad \times \int_{\partial_{k_1}^{\varepsilon_1} \cdots \partial_{k_j}^{\varepsilon_j} \underline{K}_n} \sum_{\sigma \in \mathfrak{S}_d} h_{\sigma(1)} \left( \frac{k_1}{n} \right) \cdots h_{\sigma(d)} \left( \frac{k_d}{n} \right) 
         \tilde{\varphi}_n(\underline{x}) (\Xi_{n,a,b}^n)^{-1} \\ &\quad \times \frac{\partial^{d-j} q_n^{0,n}}{\partial x_{k_{d}} \cdots \partial x_{k_{j+1}}}(\underline{x}) \ \mathrm{d} \underline{x}_{\left[ 0, \frac{k_{1}}{n} \right]}^{a,f^{\varepsilon_{1}}\left( \frac{k_{1}}{n} \right)} \prod_{i=2}^{j} \ \mathrm{d} \underline{x}_{\left[ \frac{k_{i-1}}{n}, \frac{k_{i}}{n} \right]}^{f^{\varepsilon_{i-1}}\left( \frac{k_{i-1}}{n} \right),f^{\varepsilon_i}\left( \frac{k_{i}}{n} \right)} \ \mathrm{d} \underline{x}_{\left[ \frac{k_{j}}{n}, 1 \right]}^{f^{\varepsilon_j}\left( \frac{k_{j}}{n} \right), b}.
    \end{aligned}
\end{equation}
In addition, From the assumption (\ref{phi_cond}) for $\varphi$ and H\"{o}lder's inequality, we establish
\begin{equation}\nonumber
    \left|\frac{\partial^d \tilde{\varphi}_n}{\partial x_{k_1} \cdots \partial x_{k_d}} (\underline{x})\right| \lesssim \frac{1}{n^{d-1+\delta}}~\text{for some}~\delta > 0,
\end{equation}
so it should be noted that when performing integration by parts, there is no need to consider the calculation of pseudo-diagonal terms.
The pseudo-diagonal terms refer to the $d$-tuple $(k_1, \dots, k_d)$ such that $\exists i \neq j ~\text{s.t.}~|k_i - k_j| \leq 1$.

When summarizing the results of the boundary terms calculation, we employ a technique to simplify the expressions.
To achieve conciseness, we interchange the formal indices $k_1, k_2, \dots, k_d$ within each term, such that $k_1 < k_2 < \dots < k_j$ and $k_{j+1}, k_{j+2}, \dots, k_d$ remain in an arbitrary order.
Subsequently, the terms can be combined as $\int \sum h_{\ell_1}(k_1/n) h_{\ell_2}(k_2/n) \dots h_{\ell_d}(k_d/n) \times \Theta$.
Here, $\Theta$ depends on the arrangement of $k_1, k_2, \dots, k_d$, but it is independent of the order of $k_{j+1}, k_{j+2}, \dots, k_d$.
In this context, the ranges of $\ell_1, \ell_2, \dots, \ell_d$ within the set $\Sigma$ can be chosen with certain degree of freedom.
The key observation is that the trade-off between the degree of freedom and that the $\Theta$ is independent of the arrangement of $k_{j+1}, k_{j+2}, \dots, k_d$.
In particular, the cardinality of the set of ranges within $\Sigma$ is ${}_d \mathrm{C}_j \cdot |\mathfrak{S}_j| = d! / ((d-j)!)$.
Here, ${}_d \mathrm{C}_j$ represents the binomial coefficient, denoting the number of ways to choose $j$ elements from a set of $d$ with a fixed order, and $|\mathfrak{S}_j|$ is the number of permutations of $j$ elements.
Hence, since there are $(d-j)!$ degree of freedom for arranging $k_{j+1}, k_{j+2}, \dots, k_d$, the entire expression is copied $(d-j)!$ times.
By formally reassigning indices, we can then align the ranges within $\Sigma$ with $\mathfrak{S}_d$.
This technique, akin to Symmetrization, shares similarities with the procedure used in considering the kernel representation of multiple Wiener integrals.
It is worth noting that the reason why $\Theta$ does not depend on the order of $k_{j+1}, k_{j+2}, \dots, k_d$ will be explained later in conjunction with the expression for the higher-order derivatives of $q_n^{0,n}$ and the resulting consequences.
This occurs because we introduced the notation for Symmetrization somewhat prematurely for the sake of simplicity in the expressions.

For $q_n^{0,n}(\underline{x})$,
\begin{equation}\nonumber
    \frac{\partial q_n^{0,n}}{\partial x_k}(\underline{x}) = n(x_{k+1} - 2x_k + x_{k-1}) q_n^{0,n}(\underline{x}),
\end{equation}
\begin{equation}\nonumber
    \begin{aligned}
        &\frac{\partial^2 q_n^{0,n}}{\partial x_{k_2} \partial x_{k_1}}(\underline{x}) \\
        &= nq_n^{0,n}(\underline{x}) \frac{\partial (x_{k_1 +1} - 2x_{k_1} + x_{k_1 -1})}{\partial x_{k_2}} + n(x_{k_1 +1} - 2x_{k_1} + x_{k_1 -1}) \frac{\partial q_n^{0,n}}{\partial x_{k_2}}(\underline{x}) \\
        &= nq_n^{0,n}(\underline{x}) (\delta_{k_1 +1, k_2} -2\delta_{k_1,k_2} + \delta_{k_1 -1,k_2}) + q_n^{0,n}(\underline{x}) \prod_{i=1}^2 n(x_{k_i +1} - 2x_{k_i} + x_{k_i -1}),
    \end{aligned}
\end{equation}
where $\delta_{\alpha,\beta}$ is the Kronecker delta.
In particular, when $(k_1, \dots, k_d)$ is not pseudo-diagonal, we have
\begin{equation}\nonumber
    \frac{\partial^j q_n^{0,n}}{\partial x_{k_j} \cdots \partial x_{k_1}}(\underline{x}) = q_n^{0,n}(\underline{x}) \prod_{i=1}^j n(x_{k_i +1} - 2x_{k_i} + x_{k_i -1})~\text{for}~j = 1, \dots, d.
\end{equation}
Then,
\begin{equation}\nonumber
    \begin{aligned}
        &\sum_{k_1, \dots , k_d \in \mathcal{K}_{n-1}^{[1]}} h_1 \left( \frac{k_1}{n} \right) \cdots h_d \left( \frac{k_d}{n} \right) \int_{\underline{K}_n} \frac{\partial^d \tilde{\varphi}_n}{\partial x_{k_1} \cdots \partial x_{k_d}} (\underline{x}) (\Xi_{n,a,b}^n)^{-1} q_n^{0,n}(\underline{x}) \ \mathrm{d} \underline{x}^{a,b} \\
        &= (-1)^d \sum_{k_1, \dots , k_d \in \mathcal{K}_{n-1}^{[1]}} h_1 \left( \frac{k_1}{n} \right) \cdots h_d \left( \frac{k_d}{n} \right) \\
        &\quad \times \int_{\underline{K}_n} \tilde{\varphi}_n(\underline{x}) (\Xi_{n,a,b}^n)^{-1} q_n^{0,n}(\underline{x}) \prod_{i=1}^{d} n(x_{k_i + 1} - 2x_{k_i} + x_{k_i -1}) \ \mathrm{d} \underline{x}^{a,b} \\
        &\quad + \text{BD}_n^{(1)} + \text{BD}_n^{(2)} + \dots + \text{BD}_n^{(d)},
    \end{aligned}
\end{equation}
\begin{equation}\nonumber
    \begin{aligned}
        \text{BD}_n^{(j)} &= \sum_{\varepsilon_1 , \dots , \varepsilon_j \in \{ 1, -1 \}} (-1)^{d-j} \varepsilon_1 \cdots \varepsilon_j \sum_{\substack{(k_1, \dots , k_j) \in \mathcal{K}_{n-1}^{[j]} \\ k_{j+1}, \dots , k_d \in \mathcal{K}_{n-1}^{[1]}}} \frac{1}{(d-j)!} \\ &\quad \times \int_{\partial_{k_1}^{\varepsilon_1} \cdots \partial_{k_j}^{\varepsilon_j} \underline{K}_n} \sum_{\sigma \in \mathfrak{S}_d} h_{\sigma(1)} \left( \frac{k_1}{n} \right) \cdots h_{\sigma(d)} \left( \frac{k_d}{n} \right) \\
        &\quad \times \tilde{\varphi}_n(\underline{x}) (\Xi_{n,a,b}^n)^{-1} q_n^{0,n}(\underline{x}) \prod_{i=j+1}^d n(x_{k_i + 1} - 2x_{k_i} + x_{k_i -1}) \\
        &\quad \times \ \mathrm{d} \underline{x}_{\left[ 0, \frac{k_{1}}{n} \right]}^{a,f^{\varepsilon_{1}}\left( \frac{k_{1}}{n} \right)} \prod_{i=2}^{j} \ \mathrm{d} \underline{x}_{\left[ \frac{k_{i-1}}{n}, \frac{k_{i}}{n} \right]}^{f^{\varepsilon_{i-1}}\left( \frac{k_{i-1}}{n} \right),f^{\varepsilon_i}\left( \frac{k_{i}}{n} \right)} \ \mathrm{d} \underline{x}_{\left[ \frac{k_{j}}{n}, 1 \right]}^{f^{\varepsilon_j}\left( \frac{k_{j}}{n} \right), b}.
    \end{aligned}
\end{equation}
In addition, for $\underline{x} = (x_j)_{j=k_1}^{k_2} \in \mathbb{R}^{k_2 - k_1 +1}$, we have
\begin{equation}\nonumber
    P_{\left[ \frac{k_1}{n}, \frac{k_2}{n} \right]}^{a,b} \circ \pi_{n,1}^{-1}(\ \mathrm{d} \underline{x}) = (\Xi_{n,a,b}^{k_2 - k_1})^{-1} q_n^{k_1,k_2}(\underline{x}) \delta_a(\ \mathrm{d} x_{k_1}) \prod_{i = k_1 +1}^{k_2 -1} \ \mathrm{d} x_i \delta_b(\ \mathrm{d} x_{k_2})
\end{equation}
Here, $P_{[t_1,t_2]}^{a,b}$ stands for the distribution of the one-dimensional Brownian bridge from $a$ to $b$ on $[t_1,t_2]$.
Therefore, we can express
\begin{equation}\nonumber
    \begin{aligned}
        &\sum_{k_1, \dots , k_d \in \mathcal{K}_{n-1}^{[1]}} h_1 \left( \frac{k_1}{n} \right) \cdots h_d \left( \frac{k_d}{n} \right) \int_{\underline{K}_n} \frac{\partial^d \tilde{\varphi}_n}{\partial x_{k_1} \cdots \partial x_{k_d}} (\underline{x}) (\Xi_{n,a,b}^n)^{-1} q_n^{0,n}(\underline{x}) \ \mathrm{d} \underline{x}^{a,b} \\
        &= (-1)^d \sum_{k_1, \dots , k_d \in \mathcal{K}_{n-1}^{[1]}} h_1 \left( \frac{k_1}{n} \right) \cdots h_d \left( \frac{k_d}{n} \right) \\
        &\quad \times \int_{\underline{K}_n} \tilde{\varphi}_n(\underline{x}) \prod_{i=1}^{d} n(x_{k_i + 1} - 2x_{k_i} + x_{k_i -1}) P^{a,b} \circ \pi_{n,1}^{-1}(\ \mathrm{d} \underline{x}) \\
        &\quad + \text{BD}_n^{(1)} + \text{BD}_n^{(2)} + \dots + \text{BD}_n^{(d)},
    \end{aligned}
\end{equation}
\begin{equation}\nonumber
    \begin{aligned}
        &\text{BD}_n^{(j)} \\
        &= \sum_{\varepsilon_1 , \dots , \varepsilon_j \in \{ 1, -1 \}} (-1)^{d-j} \varepsilon_1 \cdots \varepsilon_j \sum_{\substack{(k_1, \dots , k_j) \in \mathcal{K}_{n-1}^{[j]} \\ k_{j+1}, \dots , k_d \in \mathcal{K}_{n-1}^{[1]}}} \frac{1}{(d-j)!} \\
        &\quad \times \int_{\partial_{k_1}^{\varepsilon_1} \cdots \partial_{k_j}^{\varepsilon_j} \underline{K}_n} \sum_{\sigma \in \mathfrak{S}_d} h_{\sigma(1)} \left( \frac{k_1}{n} \right) \cdots h_{\sigma(d)} \left( \frac{k_d}{n} \right) \\
        &\quad \times \tilde{\varphi}_n(\underline{x}) (\Xi_{n,a,b}^n)^{-1} \Xi_{n,a,f^{\varepsilon_1}\left( \frac{k_1}{n} \right)}^{k_1} \left( \prod_{\ell=2}^j \Xi_{n,f^{\varepsilon_{\ell -1}}\left( \frac{k_{\ell -1}}{n} \right), f^{\varepsilon_{\ell}}\left( \frac{k_{\ell}}{n} \right)}^{k_{\ell} - k_{\ell -1}} \right) \Xi_{n,f^{\varepsilon_j}\left( \frac{k_j}{n} \right), b}^{n-k_j} \\
        &\quad \times \prod_{i=j+1}^d n(x_{k_i + 1} - 2x_{k_i} + x_{k_i -1}) \\
        &\quad \times P_{\left[ 0, \frac{k_1}{n} \right]}^{a, f^{\varepsilon_1} \left( \frac{k_1}{n} \right)} \circ \pi_{n,1}^{-1}(\ \mathrm{d} \underline{x}_{\left[ 0, \frac{k_1}{n} \right]}) \\
        &\quad \times \left( \prod_{i=2}^{j} P_{\left[ \frac{k_{i-1}}{n}, \frac{k_i}{n} \right]}^{f^{\varepsilon_{i-1}} \left( \frac{k_{i-1}}{n} \right), f^{\varepsilon_{i}} \left( \frac{k_{i}}{n} \right)} \circ \pi_{n,1}^{-1}(\ \mathrm{d} \underline{x}_{\left[ \frac{k_{i-1}}{n}, \frac{k_i}{n} \right]}) \right) \\
        &\quad \times P_{\left[ \frac{k_j}{n}, 1 \right]}^{f^{\varepsilon_j} \left( \frac{k_j}{n} \right), b} \circ \pi_{n,1}^{-1}(\ \mathrm{d} \underline{x}_{\left[ \frac{k_j}{n}, 1 \right]}).
    \end{aligned}
\end{equation}
For each term in this finite-dimensional approximation, we consider taking the limit as $k_1/n \to t_1, k_2/n \to t_2, \dots, k_d/n \to t_d$, and $n \to \infty$.
Basically, the same limit operation as Funaki--Ishitani~\cite{funaki_ishitani} should be applied.
Therefore, we will discuss the terms that appear newly for the case of $d (\geq 2)$-th order derivatives.
Specifically, for $2 \leq i \leq j$, we evaluate
\begin{equation}\nonumber
    P_{\left[ \frac{k_{i-1}}{n}, \frac{k_i}{n} \right]}^{f^{\varepsilon_{i-1}} \left( \frac{k_{i-1}}{n} \right), f^{\varepsilon_{i}} \left( \frac{k_{i}}{n} \right)} \circ \pi_{n,1}^{-1}(\ \mathrm{d} \underline{x}_{\left[ \frac{k_{i-1}}{n}, \frac{k_i}{n} \right]}).
\end{equation}
Now,
\begin{equation}\nonumber
    \begin{aligned}
        &P_{\left[ \frac{k_{i-1}}{n}, \frac{k_i}{n} \right]}^{f^{\varepsilon_{i-1}} \left( \frac{k_{i-1}}{n} \right), f^{\varepsilon_{i}} \left( \frac{k_{i}}{n} \right)} \circ \pi_{n,1}^{-1}(\ \mathrm{d} \underline{x}_{\left[ \frac{k_{i-1}}{n}, \frac{k_i}{n} \right]}) \\
        &= P_{\left[ \frac{k_{i-1}}{n}, \frac{k_i}{n} \right]}^{f^{\varepsilon_{i-1}} \left( \frac{k_{i-1}}{n} \right), f^{\varepsilon_{i}} \left( \frac{k_{i}}{n} \right)}(\pi_{n,1}^{-1}(\underline{K}_n^{\frac{k_{i-1}}{n}, \frac{k_i}{n}})) \\
        &\quad \times P_{\left[ \frac{k_{i-1}}{n}, \frac{k_i}{n} \right]}^{f^{\varepsilon_{i-1}} \left( \frac{k_{i-1}}{n} \right), f^{\varepsilon_{i}} \left( \frac{k_{i}}{n} \right)} \circ \pi_{n,1}^{-1}(\ \mathrm{d} \underline{x}_{\left[ \frac{k_{i-1}}{n}, \frac{k_i}{n} \right]} \vert \underline{K}_n^{\frac{k_{i-1}}{n}, \frac{k_i}{n}})
    \end{aligned}
\end{equation}
can be decomposed.
The conditional probability on the right-hand side converges weakly to Brownian excursion when $\varepsilon_{i-1} = \varepsilon_i$ and to Brownian house-moving when $\varepsilon_{i-1} = - \varepsilon_i$~\cite{yanashima}.
Furthermore, due to the lemma \ref{infinitesimal_lemma}, we have
\begin{equation}\nonumber
    \lim_{n \to \infty} n P_{\left[ \frac{k_{i-1}}{n}, \frac{k_i}{n} \right]}^{f^{\varepsilon_{i-1}} \left( \frac{k_{i-1}}{n} \right), f^{\varepsilon_{i}} \left( \frac{k_{i}}{n} \right)}(\pi_{n,1}^{-1}(\underline{K}_n^{\frac{k_{i-1}}{n}, \frac{k_i}{n}})) = \frac{\Delta^2 P_{[t_{i-1},t_i]}^{f^{\varepsilon_{i-1}}(t_{i-1}),f^{\varepsilon_i}(t_i)} (K_{[t_{i-1},t_i]})}{(\sqrt{t_{i} - t_{i-1}})^2}.
\end{equation}
Also,
\begin{equation}\nonumber
    \begin{aligned}
        & (\Xi_{n,a,b}^n)^{-1} \Xi_{n,a,f^{\varepsilon_1}\left( \frac{k_1}{n} \right)}^{k_1} \left( \prod_{\ell=2}^j \Xi_{n,f^{\varepsilon_{\ell -1}}\left( \frac{k_{\ell -1}}{n} \right), f^{\varepsilon_{\ell}}\left( \frac{k_{\ell}}{n} \right)}^{k_{\ell} - k_{\ell -1}} \right) \Xi_{n,f^{\varepsilon_j}\left( \frac{k_j}{n} \right), b}^{n-k_j} \\
        &\to p_{a,b}^{(j)}(t_1,t_2, \dots, t_j; f^{\varepsilon_1}(t_1), f^{\varepsilon_2}(t_2), \dots, f^{\varepsilon_j}(t_j)) , \\
        &\quad (n \to \infty,  k_1/n \to t_1 , k_2/n \to t_2, \dots, k_j/n \to t_j).
    \end{aligned}
\end{equation}
Therefore, we obtain
\begin{equation}\nonumber
    \begin{aligned}
        &\int_{\Omega} 1_K(X^{a,b}) \nabla^{(d)} \varphi(h_1,h_2, \dots , h_d)(X^{a,b}) \ \mathrm{d} P \\
        &= (-1)^d \int_{\Omega} 1_K(X^{a,b}) \varphi(X^{a,b}) \prod_{i=1}^{d} \int_{\mathcal{T}^{[1]}} h_i''(t_i) X^{a,b}(t_i) \ \mathrm{d} t_i \ \mathrm{d} P \\
        &\quad + \text{BD}_{\infty}^{(1)} + \text{BD}_{\infty}^{(2)} + \dots + \text{BD}_{\infty}^{(d)}.
    \end{aligned}
\end{equation}
Here, for $1 \leq j \leq d$, we have
\begin{equation}\nonumber
    \begin{aligned}
        \text{BD}_{\infty}^{(j)} =& ~ (-1)^{d-j} \sum_{\varepsilon_1 , \dots , \varepsilon_j \in \{ 1, -1 \} } \varepsilon_1 \cdots \varepsilon_j \int_{\mathcal{T}^{[j]}} \frac{1}{(d-j)!} \nu_{(j)}^{\varepsilon_1 , \dots , \varepsilon_j}(t_1, \dots , t_j) \\
        &\quad \times \int_{\Omega} \varphi(Y_{t_1, \dots , t_j}^{\varepsilon_1, \dots, \varepsilon_j, (f^-,f^+)}) \sum_{\sigma \in \mathfrak{S}_d} h_{\sigma(1)}(t_1) \cdots h_{\sigma(j)}(t_j) \\
        &\quad \times \left( \prod_{i=j+1}^{d} \int_{\mathcal{T}^{[1]}} h_{\sigma(i)}''(t_{i}) Y_{t_1, \dots , t_j}^{\varepsilon_1, \dots, \varepsilon_j, (f^-,f^+)} (t_{i}) \ \mathrm{d} t_i \right) ~ \ \mathrm{d} P \ \mathrm{d} t_1 \cdots \ \mathrm{d} t_j .
    \end{aligned}
\end{equation}
It results from Ito's lemma and the compact support of $h_1, h_2, \dots, h_d$ in $(0,1)$ that
\begin{equation}\nonumber
    \begin{aligned}
        &\int_{\mathcal{T}^{[1]}} h_{\sigma(i)}''(t_{i}) Y_{t_1, \dots , t_j}^{\varepsilon_1, \dots, \varepsilon_j, (f^-,f^+)} (t_{i}) \ \mathrm{d} t_i \\
        &= - \int_{\mathcal{T}^{[1]}} h_{\sigma(i)}'(t_{i}) \ \mathrm{d} Y_{t_1, \dots , t_j}^{\varepsilon_1, \dots, \varepsilon_j, (f^-,f^+)} (t_{i}) ,
    \end{aligned}
\end{equation}
and
\begin{equation}\nonumber
    \int_{\mathcal{T}^{[1]}} h_i''(t_i) X^{a,b}(t_i) \ \mathrm{d} t_i = - \int_{\mathcal{T}^{[1]}} h_i'(t_i) \ \mathrm{d} X^{a,b}(t_i).
\end{equation}
Then, $\text{BD}^{(j)} = \text{BD}_{\infty}^{(j)}$, so the proof of Theorem \ref{main1} is concluded.

\section{Infinitesimal probabilities}
In this section, we formulate the infinitesimal probabilities used in the proof of Theorem \ref{main1}.
\begin{df}\label{infinitesimal}
    For $(a,b) \in (f^-(t_1), f^+(t_1)) \times \{ f^-(t_2), f^+(t_2) \}$ or $(a,b) \in \{f^-(t_1), f^+(t_1)\} \times ( f^-(t_2), f^+(t_2) )$, we define the first-order infinitesimal probability $\Delta P_{[t_1,t_2]}^{a,b}(K_{[t_1,t_2]})$ as follows:
    \begin{equation}\nonumber
        \begin{aligned}
            \Delta P_{[t_1,t_2]}^{a,b}(K_{[t_1,t_2]}) &= \lim_{\substack{n \to \infty \\ \frac{k_1}{n} \to t_1, \frac{k_2}{n} \to t_2}} \sqrt{k_2 - k_1} P_{\left[ \frac{k_1}{n}, \frac{k_2}{n} \right]}^{a,b} (\pi_{n,1}^{-1} (\underline{K}_n^{a,b})) \\
            &= \sqrt{t_2 - t_1} \lim_{\substack{n \to \infty \\ \frac{k_1}{n} \to t_1, \frac{k_2}{n} \to t_2}} \sqrt{n} P_{\left[ \frac{k_1}{n}, \frac{k_2}{n} \right]}^{a,b} (\pi_{n,1}^{-1} (\underline{K}_n^{a,b})).
        \end{aligned}
    \end{equation}
    For $(a,b) \in \{f^-(t_1), f^+(t_1)\} \times \{ f^-(t_2), f^+(t_2) \}$, we define the second-order infinitesimal probability $\Delta^2 P_{[t_1,t_2]}^{a,b} (K_{[t_1,t_2]})$ as follows:
    \begin{equation}\nonumber
        \begin{aligned}
            \Delta^2 P_{[t_1,t_2]}^{a,b}(K_{[t_1,t_2]}) &= \lim_{\substack{n \to \infty \\ \frac{k_1}{n} \to t_1, \frac{k_2}{n} \to t_2}} (\sqrt{k_2 - k_1})^2 P_{\left[ \frac{k_1}{n}, \frac{k_2}{n} \right]}^{a,b} (\pi_{n,1}^{-1} (\underline{K}_n^{a,b})) \\
            &= (\sqrt{t_2 - t_1})^2 \lim_{\substack{n \to \infty \\ \frac{k_1}{n} \to t_1, \frac{k_2}{n} \to t_2}} (\sqrt{n})^2 P_{\left[ \frac{k_1}{n}, \frac{k_2}{n} \right]}^{a,b} (\pi_{n,1}^{-1} (\underline{K}_n^{a,b})).
        \end{aligned}
    \end{equation}
    For $a \in \{f^-(t), f^+(t)\}$, we define the first-order infinitesimal probability $\Delta P_{[t,1]}^{a}(K_{[t,1]})$ as follows:
    \begin{equation}\nonumber
        \begin{aligned}
            \Delta P_{[t,1]}^{a}(K_{[t,1]}) &= \lim_{\substack{n \to \infty \\ \frac{k}{n} \to t}} \sqrt{n - k} P_{\left[ \frac{k}{n}, 1 \right]}^{a} (\pi_{n,1}^{-1} (\underline{K}_n^{a})) \\
            &= \sqrt{1 - t} \lim_{\substack{n \to \infty \\ \frac{k}{n} \to t}} \sqrt{n} P_{\left[ \frac{k}{n}, 1 \right]}^{a} (\pi_{n,1}^{-1} (\underline{K}_n^{a})).
        \end{aligned}
    \end{equation}
\end{df}

\begin{lem}\label{infinitesimal_lemma}
    The limits in Definition \ref{infinitesimal} converge, and their limits are as follows.
    For $a \in (f^-(t_1), f^+(t_1))$ and $b = f^{\pm}(t_2)$:
    \begin{equation}\nonumber
        \begin{aligned}
            &\Delta P_{[t_1,t_2]}^{a,b}(K_{[t_1,t_2]}) \\
            &= \frac{\sqrt{t_2 - t_1}}{\sqrt{2}} ~ \overline{c}_{|a - f^{\pm}(t_1)|,0;\pm}^{t_1,t_2}(f^{\pm}) P(X_{[t_1,t_2]}^{a,f^{\pm}(t_2)} \mid_{K_{[t_1,t_2]}^{\mp}(f^{\pm})} \in K_{[t_1,t_2]}^{\pm}(f^{\mp})) .
        \end{aligned}
    \end{equation}
    For $a = f^{\pm}(t_1)$ and $b \in (f^-(t_2), f^+(t_2))$:
    \begin{equation}\nonumber
        \begin{aligned}
            &\Delta P_{[t_1,t_2]}^{a,b}(K_{[t_1,t_2]}) \\
            &= \frac{\sqrt{t_2 - t_1}}{\sqrt{2}} ~ \overline{c}_{0, |b - f^{\pm}(t_2)|;\pm}^{t_1,t_2}(f^{\pm}) P(X_{[t_1,t_2]}^{f^{\pm}(t_1),b} \mid_{K_{[t_1,t_2]}^{\mp}(f^{\pm})} \in K_{[t_1,t_2]}^{\pm}(f^{\mp})) .
        \end{aligned}
    \end{equation}
    For $a = f^{\varepsilon_1}(t_1)$, $b = f^{\varepsilon_2}(t_2)$, and any $\tau \in (t_1, t_2)$:
    \begin{equation}\nonumber
        \begin{aligned}
            &\Delta^2 P_{[t_1,t_2]}^{a,b}(K_{[t_1,t_2]}) \\
            =&~ \frac{(\sqrt{t_2 - t_1})^2}{2} \\
            &~ \times \int_{f^-(\tau)}^{f^+(\tau)} \overline{c}_{0,|\alpha - f^{\varepsilon_1}(\tau)|;\varepsilon_1}^{t_1,\tau}(f^{\varepsilon_1}) P(X_{[t_1,\tau]}^{f^{\varepsilon_1}(t_1), \alpha} \mid_{K_{[t_1,\tau]}^{-\varepsilon_1}(f^{\varepsilon_1})} \in K_{[t_1,\tau]}^{\varepsilon_1}(f^{-\varepsilon_1})) \\
            &~ \times \overline{c}_{|\alpha - f^{\varepsilon_2}(\tau)|,0;\varepsilon_2}^{\tau,t_2}(f^{\varepsilon_2}) P(X_{[\tau,t_2]}^{\alpha, f^{\varepsilon_2}(t_2)} \mid_{K_{[\tau,t_2]}^{-\varepsilon_2}(f^{\varepsilon_2})} \in K_{[\tau,t_2]}^{\varepsilon_2}(f^{-\varepsilon_2})) \\
            &~ \times P(X_{[t_1,t_2]}^{f^{\varepsilon_1}(t_1), f^{\varepsilon_2}(t_2)}(\tau) \in \ \mathrm{d} \alpha) \\
            =&~ (\sqrt{t_2 - t_1})^2 \int_{f^-(\tau)}^{f^+(\tau)} \frac{\Delta P_{[t_1,\tau]}^{f^{\varepsilon_1}(t_1), \alpha}}{\sqrt{\tau - t_1}} \frac{\Delta P_{[\tau, t_2]}^{\alpha, f^{\varepsilon_2}(t_2)}}{\sqrt{t_2 - \tau}} P(X_{[t_1,t_2]}^{f^{\varepsilon_1}(t_1), f^{\varepsilon_2}(t_2)}(\tau) \in \ \mathrm{d} \alpha) .
        \end{aligned}
    \end{equation}
    For $a = f^{\pm}(t)$:
    \begin{equation}\nonumber
        \Delta P_{[t,1]}^{a}(K_{[t,1]}) = \frac{\sqrt{1-t}}{\sqrt{2}} ~ \overline{d}_{\pm}^{t,1}(f^{\pm}) P(X_{[t,1]}^{f^{\pm}(t)} \mid_{K_{[t,1]}^{\mp}(f^{\pm})} \in K_{[t,1]}^{\pm}(f^{\mp})) .
    \end{equation}
    Here,
    \begin{equation}\nonumber
        \overline{c}_{\alpha,\beta;\pm 1}^{t_1,t_2}(f) := \frac{2|\alpha - \beta|}{t_2 - t_1} \frac{E\left[ Z_{[t_1,t_2]}^{\pm f} \left( X_{[t_1,t_2]}^{\alpha,\beta} \vert_{K_{[t_1,t_2]}^+(0)} \right) \right]}{E\left[ Z_{[t_1,t_2]}^{\pm f} \left( X_{[t_1,t_2]}^{\alpha,\beta} \right) \right]}
    \end{equation}
    and
    \begin{equation}\nonumber
        \overline{d}_{\pm 1}^{t_1,t_2}(f) := \frac{\sqrt{2}}{\sqrt{\pi (t_2 - t_1)}} E[Z_{[t_1,t_2]}^{\pm f}(X_{[t_1,t_2]}^{0} \mid_{K_{[t_1,t_2]}^{+}(0)})]
    \end{equation}
    are defined.
    Also, $K_{[t_1,t_2]}^{\pm 1} := K_{[t_1,t_2]}^{\pm}$.
\end{lem}

\subsection{Proof of Lemma \ref{infinitesimal_lemma}}
For the case of first-order infinitesimal probability, it has already been discussed by Funaki--Ishitani~\cite{funaki_ishitani}.
In the case of second-order infinitesimal probability, by fixing an arbitrary $\theta_{1,2} \in \mathcal{K}_{n-1}^{[1]}$ such that $k_1 < \theta_{1,2} < k_2$, and using the Markov property, we have
\begin{equation}\nonumber
    \begin{aligned}
        &nP_{\left[ \frac{k_1}{n}, \frac{k_2}{n} \right]}^{f^{\varepsilon_1}\left( \frac{k_1}{n} \right), f^{\varepsilon_2}\left( \frac{k_2}{n} \right)} (\pi_{n,1}^{-1}(\underline{K}_n^{\frac{k_1}{n}, \frac{k_2}{n}})) \\
        &= \int_{f^-(\frac{\theta_{1,2}}{n})}^{f^+(\frac{\theta_{1,2}}{n})} \sqrt{n} P_{\left[ \frac{k_1}{n}, \frac{\theta_{1,2}}{n} \right]}^{f^{\varepsilon_1}\left( \frac{k_1}{n} \right), \alpha } (\pi_{n,1}^{-1}(\underline{K}_n^{\frac{k_1}{n}, \frac{\theta_{1,2}}{n}})) \sqrt{n} P_{\left[ \frac{\theta_{1,2}}{n}, \frac{k_2}{n} \right]}^{\alpha , f^{\varepsilon_2}\left( \frac{k_2}{n} \right)} (\pi_{n,1}^{-1}(\underline{K}_n^{\frac{\theta_{1,2}}{n}, \frac{k_2}{n}})) \\
        &~ \times P_{\left[ \frac{k_1}{n}, \frac{k_2}{n} \right]}^{f^{\varepsilon_1}\left( \frac{k_1}{n} \right), f^{\varepsilon_2}\left( \frac{k_2}{n} \right)}(x ( \frac{\theta_{1,2}}{n} ) \in \ \mathrm{d} \alpha)
    \end{aligned}
\end{equation}
The terms appearing in the integrand
\begin{equation}\nonumber
    \sqrt{n} P_{\left[ \frac{k_1}{n}, \frac{\theta_{1,2}}{n} \right]}^{f^{\varepsilon_1}\left( \frac{k_1}{n} \right), \alpha } (\pi_{n,1}^{-1}(\underline{K}_n^{\frac{k_1}{n}, \frac{\theta_{1,2}}{n}}))
\end{equation}
and
\begin{equation}\nonumber
    \sqrt{n} P_{\left[ \frac{\theta_{1,2}}{n}, \frac{k_2}{n} \right]}^{\alpha , f^{\varepsilon_2}\left( \frac{k_2}{n} \right)} (\pi_{n,1}^{-1}(\underline{K}_n^{\frac{\theta_{1,2}}{n}, \frac{k_2}{n}}))
\end{equation}
can be expressed using the formula for the first-order infinitesimal probability by taking the limit $k_1/n \to t_1, k_2/n \to t_2, \theta_{1,2}/n \to \tau$, and $n \to \infty$.

\section{Proof of Theorem \ref{main2}}
Fundamentally, the proof can be conducted by the same manipulations as Theorem \ref{main1}.
First, we consider the integral on $\Omega$ as an integral on the continuous function space $C([0,1])$ and approximate it with integrals over finite-dimensional spaces using Proposition \ref{polygonal}.
\begin{equation}\nonumber
    \begin{aligned}
        &\int_{\Omega} 1_K(X^{a}) \nabla^{(d)} \varphi(h_1,h_2, \dots , h_d)(X^{a}) \ \mathrm{d} P \\
        &= \lim_{n \to \infty} \int_{\underline{K}_n} \sum_{k_1, \dots , k_d \in \mathcal{K}_{n-1}^{[1]}} \frac{\partial^d \tilde{\varphi}_n}{\partial x_{k_1} \cdots \partial x_{k_d}} (\underline{x}) (\pi_{n,1} h_1)_{k_1} \cdots (\pi_{n,1} h_d)_{k_d} P^{a} \circ \pi_{n,1}^{-1}(\ \mathrm{d} \underline{x})
    \end{aligned}
\end{equation}
Now,
\begin{equation}\nonumber
    P^a \circ \pi_{n,1}^{-1}(\ \mathrm{d} \underline{x}) = \left( \frac{n}{2\pi} \right)^{\frac{n}{2}} q_n^{0,n}(\underline{x}) \ \mathrm{d} \underline{x}^a
\end{equation}
is used, and by applying the integration by parts formula on $\mathbb{R}^{n-1}$, we obtain
\begin{equation}\nonumber
    \begin{aligned}
        &\int_{\underline{K}_n} \sum_{k_1, \dots , k_d \in \mathcal{K}_{n-1}^{[1]}} \frac{\partial^d \tilde{\varphi}_n}{\partial x_{k_1} \cdots \partial x_{k_d}} (\underline{x}) (\pi_{n,1} h_1)_{k_1} \cdots (\pi_{n,1} h_d)_{k_d} P^{a} \circ \pi_{n,1}^{-1}(\ \mathrm{d} \underline{x}) \\
        &= \sum_{k_1, \dots , k_d \in \mathcal{K}_{n-1}^{[1]}} h_1 \left( \frac{k_1}{n} \right) \cdots h_d \left( \frac{k_d}{n} \right) \\
        &\quad \times \int_{\underline{K}_n} \tilde{\varphi}_n(\underline{x}) \prod_{i=1}^{d} n(x_{k_i + 1} - 2x_{k_i} + x_{k_i -1}) P^{a} \circ \pi_{n,1}^{-1}(\ \mathrm{d} \underline{x}) \\
        &\quad + \text{BD}_n^{(1)} + \text{BD}_n^{(2)} + \dots + \text{BD}_n^{(d)} .
    \end{aligned}
\end{equation}
Here,
\begin{equation}\nonumber
    \begin{aligned}
        &\text{BD}_n^{(j)} \\
        &= \sum_{\varepsilon_1 , \dots , \varepsilon_j \in \{ 1, -1 \}} (-1)^{d-j} \varepsilon_1 \cdots \varepsilon_j \sum_{\substack{(k_1, \dots , k_j) \in \mathcal{K}_{n-1}^{[j]} \\ k_{j+1}, \dots , k_d \in \mathcal{K}_{n-1}^{[1]}}} \frac{1}{(d-j)!} \\ &\quad \times \int_{\partial_{k_1}^{\varepsilon_1} \cdots \partial_{k_j}^{\varepsilon_j} \underline{K}_n} \sum_{\sigma \in \mathfrak{S}_d} h_{\sigma(1)} \left( \frac{k_1}{n} \right) \cdots h_{\sigma(d)} \left( \frac{k_d}{n} \right) \\
        &\quad \times (\Xi_{n,a,f^{\varepsilon_1}\left( \frac{k_1}{n} \right)}^{k_1}) \left( \prod_{i=2}^j \Xi_{n,f^{\varepsilon_{i-1}}\left( \frac{k_{i-1}}{n} \right), f^{\varepsilon_i}\left( \frac{k_i}{n} \right)}^{k_i - k_{i-1}} \right) \left( \frac{n}{2\pi} \right)^{\frac{k_j}{2}} \\
        &\quad \times \left( \prod_{i=j+1}^{d} n (x(k_{i}+1) - 2x(k_{i}) + x(k_{i}-1)) \right) \tilde{\varphi}_n(\underline{x}) \\
        &\quad \times P_{\left[ 0, \frac{k_1}{n} \right]}^{a, f^{\varepsilon_1} \left( \frac{k_1}{n} \right)} \circ \pi_{n,1}^{-1}(\ \mathrm{d} \underline{x}_{\left[ 0, \frac{k_1}{n} \right]}) \left( \prod_{i=2}^{j} P_{\left[ \frac{k_{i-1}}{n}, \frac{k_i}{n} \right]}^{f^{\varepsilon_{i-1}} \left( \frac{k_{i-1}}{n} \right), f^{\varepsilon_{i}} \left( \frac{k_{i}}{n} \right)} \circ \pi_{n,1}^{-1}(\ \mathrm{d} \underline{x}_{\left[ \frac{k_{i-1}}{n}, \frac{k_i}{n} \right]}) \right) \\ &\quad \times P_{\left[ \frac{k_j}{n}, 1 \right]}^{f^{\varepsilon_j} \left( \frac{k_j}{n} \right)} \circ \pi_{n,1}^{-1}(\ \mathrm{d} \underline{x}_{\left[ \frac{k_j}{n}, 1 \right]}) .
    \end{aligned}
\end{equation}
Now, for $\underline{x} = (x_j)_{j=k_1}^{k_2} \in \mathbb{R}^{k_2 - k_1 +1}$,
\begin{equation}\nonumber
    P_{\left[ \frac{k_1}{n}, \frac{k_2}{n} \right]}^{a} \circ \pi_{n,1}^{-1}(\ \mathrm{d} \underline{x}) = \left( \frac{n}{2\pi} \right)^{\frac{k_2 - k_1}{2}} q_n^{k_1,k_2}(\underline{x}) \delta_a(\ \mathrm{d} x_{k_1}) \prod_{i = k_1 +1}^{k_2 -1} \ \mathrm{d} x_i \delta_b(\ \mathrm{d} x_{k_2})
\end{equation}
is used, and hence,
\begin{equation}\nonumber
    \begin{aligned}
        &(\Xi_{n,a,f^{\varepsilon_1}\left( \frac{k_1}{n} \right)}^{k_1}) \left( \prod_{i=2}^j \Xi_{n,f^{\varepsilon_{i-1}}\left( \frac{k_{i-1}}{n} \right), f^{\varepsilon_i}\left( \frac{k_i}{n} \right)}^{k_i - k_{i-1}} \right) \left( \frac{n}{2\pi} \right)^{\frac{k_j}{2}} \\
        &\to p_a^{(j)}(t_1, t_2, \dots, t_j; f^{\varepsilon_1}(t_1), f^{\varepsilon_2}(t_2), \dots, f^{\varepsilon_j}(t_j)) ,\\
        &\quad (n \to \infty,  k_1/n \to t_1 , k_2/n \to t_2, \dots, k_j/n \to t_j) .
    \end{aligned}
\end{equation}
Therefore, by taking limits as in the proof of Theorem \ref{main1}, the proof is concluded.

\begin{appendices}

\section{Brownian house-moving}
A stochastic process called Brownian house-moving can be constructed as a weak convergence limit of a Brownian bridge conditioned to stay between two curves~\cites{ishitani1, ishitani_hatakenaka_suzuki}.
\begin{prop}\label{house-moving1}
    Assume $f^-,\ f^+ \in C^2([t_1,t_2], \mathbb{R})$.
    Then, there exists a $\mathbb{R}$-valued continuous Markov process $H^{f^- \to f^+} = \{ H^{f^- \to f^+}(t) \}_{t \in [t_1,t_2]}$, such that for any $\mathbb{R}$-valued bounded continuous function $F$ on $C([t_1,t_2], \mathbb{R})$, the following convergence holds:
    \begin{equation}\nonumber
        E[F(H^{f^- \to f^+})] = \lim_{\varepsilon \downarrow 0} E[F(B_{[t_1,t_2]}^{a,b} \vert_{K_{[t_1,t_2]}(f^- - \eta^{-}(\varepsilon), f^+ + \eta^{+}(\varepsilon))})]
    \end{equation}
    where we denote $H^{f^- \to f^+}$ as $X^{a,b,(f^-,f^+)}$.
    In particular, it is referred to as Brownian excursion when $a = f^{\pm}(t_1)$ and $b = f^{\pm}(t_2)$, and as Brownian house-moving when $a = f^{\pm}(t_1)$ and $b = f^{\mp}(t_2)$.
\end{prop}
The probability density function and transition density of Brownian house-moving are already known~\cite{ishitani_hatakenaka_suzuki}.
Additionally, the construction of Brownian house-moving through random walk approximation is also known~\cite{yanashima}.
\begin{prop}
    Assume $f^-,\ f^+ \in C^2([t_1,t_2], \mathbb{R})$.
    Then, there exists a $\mathbb{R}$-valued continuous Markov process $H^{f^- \to f^+} = \{ H^{f^- \to f^+}(t) \}_{t \in [t_1,t_2]}$, such that for any $\mathbb{R}$-valued bounded continuous function $F$ on $C([t_1,t_2], \mathbb{R})$, the following convergence holds:
    \begin{equation}\nonumber
        \lim_{n \to \infty} E[F(\pi_{n,1} \circ B^{a \to b}) \vert \underline{K}_n(f^-, f^+)] = E[F(H^{f^- \to f^+})] .
    \end{equation}
    The weak limit $H^{f^- \to f^+}$ obtained here has the same probability law as the one obtained in Proposition \ref{house-moving1}.
\end{prop}

\section{Approximations by smooth curves}\label{smooth_curves}
According to Lemma 8.1 in~\cite{funaki_ishitani}, if $f_n^{\pm} \to f^{\pm}$ in $W^{1,p}([0,1])$ for $p > 2$, then,
\begin{equation}\nonumber
    \Delta P_{[t_1,t_2]}^{a,b}(K_{[t_1,t_2]}(f_n^-, f_n^+)) \to \Delta P_{[t_1,t_2]}^{a,b}(K_{[t_1,t_2]}(f^-, f^+)) \quad (n \to \infty) .
\end{equation}
Furthermore, based on Lemma \ref{infinitesimal_lemma} in this study,
\begin{equation}\nonumber
    \begin{aligned}
        &\Delta^2 P_{[t_1,t_2]}^{a,b}(K_{[t_1,t_2]}) \\
        &= (\sqrt{t_2 - t_1})^2 \int_{f^-(\tau)}^{f^+(\tau)} \frac{\Delta P_{[t_1,\tau]}^{f^{\varepsilon_1}(t_1), \alpha}}{\sqrt{\tau - t_1}} \frac{\Delta P_{[\tau, t_2]}^{\alpha, f^{\varepsilon_2}(t_2)}}{\sqrt{t_2 - \tau}} P(X_{[t_1,t_2]}^{f^{\varepsilon_1}(t_1), f^{\varepsilon_2}(t_2)}(\tau) \in \ \mathrm{d} \alpha) .
    \end{aligned}
\end{equation}
Therefore, for the second-order infinitesimal probability as well, we have
\begin{equation}\nonumber
    \Delta^2 P_{[t_1,t_2]}^{a,b}(K_{[t_1,t_2]}(f_n^-, f_n^+)) \to \Delta^2 P_{[t_1,t_2]}^{a,b}(K_{[t_1,t_2]}(f^-, f^+)) \quad (n \to \infty) .
\end{equation}

\section{Detailed calculations in the proof of Theorem \ref{main1}}
We perform detailed calculations for the case $d=2$.
When performing an approximation in a finite-dimensional space, we have
\begin{equation}\nonumber
    \begin{aligned}
        &\int_{\Omega} 1_K(X^{a,b}) \nabla^{(2)} \varphi(h_1,h_2)(X^{a,b}) \ \mathrm{d} P \\
        &= \lim_{n \to \infty} \int_{\underline{K}_n} \sum_{k_1, k_2 \in \mathcal{K}_{n-1}^{[1]}} \frac{\partial^2 \tilde{\varphi}_n}{\partial x_{k_1} \partial x_{k_2}} (\underline{x}) (\pi_{n,1} h_1)_{k_1} (\pi_{n,1} h_2)_{k_2} P^{a,b} \circ \pi_{n,1}^{-1}(\ \mathrm{d} \underline{x}),
    \end{aligned}
\end{equation}
and
\begin{equation}\nonumber
    \begin{aligned}
        &\int_{\underline{K}_n} \sum_{k_1, k_2 \in \mathcal{K}_{n-1}^{[1]}} \frac{\partial^2 \tilde{\varphi}_n}{\partial x_{k_1} \partial x_{k_2}} (\underline{x}) (\pi_{n,1} h_1)_{k_1} (\pi_{n,1} h_2)_{k_2} P^{a,b} \circ \pi_{n,1}^{-1}(\ \mathrm{d} \underline{x}) \\
        &= \sum_{k_1, k_2 \in \mathcal{K}_{n-1}^{[1]}} h_1 \left( \frac{k_1}{n} \right) h_2 \left( \frac{k_2}{n} \right) \int_{\underline{K}_n} \frac{\partial^2 \tilde{\varphi}_n}{\partial x_{k_1} \partial x_{k_2}} (\underline{x}) (\Xi_{n,a,b}^n)^{-1} q_n^{0,n}(\underline{x}) \ \mathrm{d} \underline{x}^{a,b},
    \end{aligned}
\end{equation}
so we proceed with the integration by parts on $\mathbb{R}^{n-1}$ for this equation.
\begin{equation}\nonumber
    \begin{aligned}
        &\sum_{k_1, k_2 \in \mathcal{K}_{n-1}^{[1]}} h_1 \left( \frac{k_1}{n} \right) h_2 \left( \frac{k_2}{n} \right) \int_{\underline{K}_n} \frac{\partial^2 \tilde{\varphi}_n}{\partial x_{k_1} \partial x_{k_2}} (\underline{x}) (\Xi_{n,a,b}^n)^{-1} q_n^{0,n}(\underline{x}) \ \mathrm{d} \underline{x}^{a,b} \\
        &= - \sum_{k_1, k_2 \in \mathcal{K}_{n-1}^{[1]}} h_1 \left( \frac{k_1}{n} \right) h_2 \left( \frac{k_2}{n} \right) \int_{\underline{K}_n} \frac{\partial \tilde{\varphi}_n}{\partial x_{k_2}} (\underline{x}) (\Xi_{n,a,b}^n)^{-1} \frac{\partial q_n^{0,n}}{\partial x_{k_1}}(\underline{x}) \ \mathrm{d} \underline{x}^{a,b} \\
        &\quad + \sum_{k_1, k_2 \in \mathcal{K}_{n-1}^{[1]}} h_1 \left( \frac{k_1}{n} \right) h_2 \left( \frac{k_2}{n} \right) \int_{\partial_{k_1}^+ \underline{K}_n} \frac{\partial \tilde{\varphi}_n}{\partial x_{k_2}} (\underline{x}) (\Xi_{n,a,b}^n)^{-1} q_n^{0,n}(\underline{x}) \\ &\qquad \qquad \qquad \qquad \qquad \qquad \qquad \qquad \times \mathrm{d} \underline{x}_{[0,k_1/n]}^{a,f^+(k_1/n)} \ \mathrm{d} \underline{x}_{[k_1/n,1]}^{f^+(k_1/n), b} \\
        &\quad - \sum_{k_1, k_2 \in \mathcal{K}_{n-1}^{[1]}} h_1 \left( \frac{k_1}{n} \right) h_2 \left( \frac{k_2}{n} \right) \int_{\partial_{k_1}^- \underline{K}_n} \frac{\partial \tilde{\varphi}_n}{\partial x_{k_2}} (\underline{x}) (\Xi_{n,a,b}^n)^{-1} q_n^{0,n}(\underline{x}) \\ &\qquad \qquad \qquad \qquad \qquad \qquad \qquad \qquad \times \mathrm{d} \underline{x}_{[0,k_1/n]}^{a,f^-(k_1/n)} \ \mathrm{d} \underline{x}_{[k_1/n,1]}^{f^-(k_1/n), b} \\
        &=: -I_n + B_n^+ - B_n^-,
    \end{aligned}
\end{equation}
where $\varphi(x) = \Phi(\langle x, \lambda_1 \rangle_H, \dots, \langle x, \lambda_{\ell} \rangle_H)$, so it suffices to consider the non-pseudo-diagonal case.
\begin{equation}\nonumber
    \begin{aligned}
        I_n &= \sum_{k_1, k_2 \in \mathcal{K}_{n-1}^{[1]}} h_1 \left( \frac{k_1}{n} \right) h_2 \left( \frac{k_2}{n} \right) \int_{\underline{K}_n} \frac{\partial \tilde{\varphi}_n}{\partial x_{k_2}} (\underline{x}) (\Xi_{n,a,b}^n)^{-1} \frac{\partial q_n^{0,n}}{\partial x_{k_1}}(\underline{x}) \ \mathrm{d} \underline{x}^{a,b} \\
        &= -\sum_{k_1, k_2 \in \mathcal{K}_{n-1}^{[1]}} h_1 \left( \frac{k_1}{n} \right) h_2 \left( \frac{k_2}{n} \right) \int_{\underline{K}_n} \tilde{\varphi}_n (\underline{x}) (\Xi_{n,a,b}^n)^{-1} \frac{\partial^2 q_n^{0,n}}{\partial x_{k_2} \partial x_{k_1}}(\underline{x}) \ \mathrm{d} \underline{x}^{a,b} \\
        &\quad + \sum_{k_1, k_2 \in \mathcal{K}_{n-1}^{[1]}} h_1 \left( \frac{k_1}{n} \right) h_2 \left( \frac{k_2}{n} \right) \int_{\partial_{k_2}^+ \underline{K}_n} \tilde{\varphi}_n (\underline{x}) (\Xi_{n,a,b}^n)^{-1} \frac{\partial q_n^{0,n}}{\partial x_{k_1}}(\underline{x}) \\ &\qquad \qquad \qquad \qquad \qquad \qquad \qquad \qquad \times \mathrm{d} \underline{x}_{[0,k_2/n]}^{a,f^+(k_2/n)} \ \mathrm{d} \underline{x}_{[k_2/n,1]}^{f^+(k_2/n), b} \\
        &\quad - \sum_{k_1, k_2 \in \mathcal{K}_{n-1}^{[1]}} h_1 \left( \frac{k_1}{n} \right) h_2 \left( \frac{k_2}{n} \right) \int_{\partial_{k_2}^- \underline{K}_n} \tilde{\varphi}_n (\underline{x}) (\Xi_{n,a,b}^n)^{-1} \frac{\partial q_n^{0,n}}{\partial x_{k_1}}(\underline{x}) \\ &\qquad \qquad \qquad \qquad \qquad \qquad \qquad \qquad \times \mathrm{d} \underline{x}_{[0,k_2/n]}^{a,f^+(k_2/n)} \ \mathrm{d} \underline{x}_{[k_2/n,1]}^{f^+(k_2/n), b}
    \end{aligned}
\end{equation}

\begin{equation}\nonumber
    \begin{aligned}
        B_n^+ &= \sum_{k_1, k_2 \in \mathcal{K}_{n-1}^{[1]}} h_1 \left( \frac{k_1}{n} \right) h_2 \left( \frac{k_2}{n} \right) \int_{\partial_{k_1}^+ \underline{K}_n} \frac{\partial \tilde{\varphi}_n}{\partial x_{k_2}} (\underline{x}) (\Xi_{n,a,b}^n)^{-1} q_n^{0,n}(\underline{x}) \\ &\qquad \qquad \qquad \qquad \qquad \qquad \qquad \qquad \times \mathrm{d} \underline{x}_{[0,k_1/n]}^{a,f^+(k_1/n)} \ \mathrm{d} \underline{x}_{[k_1/n,1]}^{f^+(k_1/n), b} \\
        &= \sum_{(k_1, k_2) \in \mathcal{K}_{n-1}^{[2]}} h_1 \left( \frac{k_1}{n} \right) h_2 \left( \frac{k_2}{n} \right) \int_{\partial_{k_1}^+ \underline{K}_n} \frac{\partial \tilde{\varphi}_n}{\partial x_{k_2}} (\underline{x}) (\Xi_{n,a,b}^n)^{-1} q_n^{0,n}(\underline{x}) \\ &\qquad \qquad \qquad \qquad \qquad \qquad \qquad \qquad \times \mathrm{d} \underline{x}_{[0,k_1/n]}^{a,f^+(k_1/n)} \ \mathrm{d} \underline{x}_{[k_1/n,1]}^{f^+(k_1/n), b} \\
        &\quad + \sum_{(k_2, k_1) \in \mathcal{K}_{n-1}^{[2]}} h_1 \left( \frac{k_1}{n} \right) h_2 \left( \frac{k_2}{n} \right) \int_{\partial_{k_1}^+ \underline{K}_n} \frac{\partial \tilde{\varphi}_n}{\partial x_{k_2}} (\underline{x}) (\Xi_{n,a,b}^n)^{-1} q_n^{0,n}(\underline{x}) \\ &\qquad \qquad \qquad \qquad \qquad \qquad \qquad \qquad \times \mathrm{d} \underline{x}_{[0,k_1/n]}^{a,f^+(k_1/n)} \ \mathrm{d} \underline{x}_{[k_1/n,1]}^{f^+(k_1/n), b} \\
        &= - \sum_{(k_1, k_2) \in \mathcal{K}_{n-1}^{[2]}} h_1 \left( \frac{k_1}{n} \right) h_2 \left( \frac{k_2}{n} \right) \int_{\partial_{k_1}^+ \underline{K}_n} \tilde{\varphi}_n (\underline{x}) (\Xi_{n,a,b}^n)^{-1} \frac{\partial q_n^{0,n}}{\partial x_{k_2}}(\underline{x}) \\ &\qquad \qquad \qquad \qquad \qquad \qquad \qquad \qquad \times \mathrm{d} \underline{x}_{[0,k_1/n]}^{a,f^+(k_1/n)} \ \mathrm{d} \underline{x}_{[k_1/n,1]}^{f^+(k_1/n), b} \\
        &\quad + \sum_{(k_1, k_2) \in \mathcal{K}_{n-1}^{[2]}} h_1 \left( \frac{k_1}{n} \right) h_2 \left( \frac{k_2}{n} \right) \int_{\partial_{k_2}^+ \partial_{k_1}^+ \underline{K}_n} \tilde{\varphi}_n (\underline{x}) (\Xi_{n,a,b}^n)^{-1} q_n^{0,n}(\underline{x}) \\
        &\qquad \qquad \qquad \qquad \times \ \mathrm{d} \underline{x}_{[0,k_1/n]}^{a,f^+(k_1/n)} \ \mathrm{d} \underline{x}_{[k_1/n, k_2/n]}^{f^+(k_1/n), f^+(k_2/n)} \ \mathrm{d} \underline{x}_{[k_2/n,1]}^{f^+(k_2/n), b} \\
        &\quad - \sum_{(k_1, k_2) \in \mathcal{K}_{n-1}^{[2]}} h_1 \left( \frac{k_1}{n} \right) h_2 \left( \frac{k_2}{n} \right) \int_{\partial_{k_2}^- \partial_{k_1}^+ \underline{K}_n} \tilde{\varphi}_n (\underline{x}) (\Xi_{n,a,b}^n)^{-1} q_n^{0,n}(\underline{x}) \\
        &\qquad \qquad \qquad \qquad \times \ \mathrm{d} \underline{x}_{[0,k_1/n]}^{a,f^+(k_1/n)} \ \mathrm{d} \underline{x}_{[k_1/n, k_2/n]}^{f^+(k_1/n), f^-(k_2/n)} \ \mathrm{d} \underline{x}_{[k_1/n,1]}^{f^-(k_2/n), b} \\
        &\quad - \sum_{(k_2, k_1) \in \mathcal{K}_{n-1}^{[2]}} h_1 \left( \frac{k_1}{n} \right) h_2 \left( \frac{k_2}{n} \right) \int_{\partial_{k_1}^+ \underline{K}_n} \tilde{\varphi}_n (\underline{x}) (\Xi_{n,a,b}^n)^{-1} \frac{\partial q_n^{0,n}}{\partial x_{k_2}}(\underline{x}) \\ &\qquad \qquad \qquad \qquad \qquad \qquad \qquad \qquad \times \mathrm{d} \underline{x}_{[0,k_1/n]}^{a,f^+(k_1/n)} \ \mathrm{d} \underline{x}_{[k_1/n,1]}^{f^+(k_1/n), b} \\
        &\quad + \sum_{(k_2, k_1) \in \mathcal{K}_{n-1}^{[2]}} h_1 \left( \frac{k_1}{n} \right) h_2 \left( \frac{k_2}{n} \right) \int_{\partial_{k_2}^+ \partial_{k_1}^+ \underline{K}_n} \tilde{\varphi}_n (\underline{x}) (\Xi_{n,a,b}^n)^{-1} q_n^{0,n}(\underline{x}) \\
        &\qquad \qquad \qquad \qquad \times \ \mathrm{d} \underline{x}_{[0,k_2/n]}^{a,f^+(k_2/n)} \ \mathrm{d} \underline{x}_{[k_2/n, k_1/n]}^{f^+(k_2/n), f^+(k_1/n)} \ \mathrm{d} \underline{x}_{[k_1/n,1]}^{f^+(k_1/n), b} \\
        &\quad - \sum_{(k_2, k_1) \in \mathcal{K}_{n-1}^{[2]}} h_1 \left( \frac{k_1}{n} \right) h_2 \left( \frac{k_2}{n} \right) \int_{\partial_{k_2}^- \partial_{k_1}^+ \underline{K}_n} \tilde{\varphi}_n (\underline{x}) (\Xi_{n,a,b}^n)^{-1} q_n^{0,n}(\underline{x}) \\
        &\qquad \qquad \qquad \qquad \times \ \mathrm{d} \underline{x}_{[0,k_2/n]}^{a,f^-(k_2/n)} \ \mathrm{d} \underline{x}_{[k_2/n, k_1/n]}^{f^-(k_2/n), f^+(k_1/n)} \ \mathrm{d} \underline{x}_{[k_1/n,1]}^{f^+(k_1/n), b} \\
    \end{aligned}
\end{equation}

\begin{equation}\nonumber
    \begin{aligned}
        &= - \sum_{k_1, k_2 \in \mathcal{K}_{n-1}^{[1]}} h_1 \left( \frac{k_1}{n} \right) h_2 \left( \frac{k_2}{n} \right) \int_{\partial_{k_1}^+ \underline{K}_n} \tilde{\varphi}_n (\underline{x}) (\Xi_{n,a,b}^n)^{-1} \frac{\partial q_n^{0,n}}{\partial x_{k_2}}(\underline{x}) \\ &\qquad \qquad \qquad \qquad \qquad \qquad \qquad \qquad \mathrm{d} \underline{x}_{[0,k_1/n]}^{a,f^+(k_1/n)} \ \mathrm{d} \underline{x}_{[k_1/n,1]}^{f^+(k_1/n), b} \\
        &\quad + \sum_{(k_1, k_2) \in \mathcal{K}_{n-1}^{[2]}} \left( h_1 \left( \frac{k_1}{n} \right) h_2 \left( \frac{k_2}{n} \right) + h_1 \left( \frac{k_2}{n} \right) h_2 \left( \frac{k_1}{n} \right) \right) \\
        &\qquad \qquad \qquad \qquad \times \int_{\partial_{k_2}^+ \partial_{k_1}^+ \underline{K}_n} \tilde{\varphi}_n (\underline{x}) (\Xi_{n,a,b}^n)^{-1} q_n^{0,n}(\underline{x}) \\
        &\qquad \qquad \qquad \qquad \times \ \mathrm{d} \underline{x}_{[0,k_1/n]}^{a,f^+(k_1/n)} \ \mathrm{d} \underline{x}_{[k_1/n, k_2/n]}^{f^+(k_1/n), f^+(k_2/n)} \ \mathrm{d} \underline{x}_{[k_2/n,1]}^{f^+(k_2/n), b} \\
        &\quad - \sum_{(k_1, k_2) \in \mathcal{K}_{n-1}^{[2]}} h_1 \left( \frac{k_1}{n} \right) h_2 \left( \frac{k_2}{n} \right) \int_{\partial_{k_2}^- \partial_{k_1}^+ \underline{K}_n} \tilde{\varphi}_n (\underline{x}) (\Xi_{n,a,b}^n)^{-1} q_n^{0,n}(\underline{x}) \\
        &\qquad \qquad \qquad \qquad \times \ \mathrm{d} \underline{x}_{[0,k_1/n]}^{a,f^+(k_1/n)} \ \mathrm{d} \underline{x}_{[k_1/n, k_2/n]}^{f^+(k_1/n), f^-(k_2/n)} \ \mathrm{d} \underline{x}_{[k_1/n,1]}^{f^-(k_2/n), b} \\
        &\quad - \sum_{(k_2, k_1) \in \mathcal{K}_{n-1}^{[2]}} h_1 \left( \frac{k_1}{n} \right) h_2 \left( \frac{k_2}{n} \right) \int_{\partial_{k_2}^- \partial_{k_1}^+ \underline{K}_n} \tilde{\varphi}_n (\underline{x}) (\Xi_{n,a,b}^n)^{-1} q_n^{0,n}(\underline{x}) \\
        &\qquad \qquad \qquad \qquad \times \ \mathrm{d} \underline{x}_{[0,k_2/n]}^{a,f^-(k_2/n)} \ \mathrm{d} \underline{x}_{[k_2/n, k_1/n]}^{f^-(k_2/n), f^+(k_1/n)} \ \mathrm{d} \underline{x}_{[k_1/n,1]}^{f^+(k_1/n), b}
    \end{aligned}
\end{equation}

\begin{equation}\nonumber
    \begin{aligned}
        B_n^- &= \sum_{k_1, k_2 \in \mathcal{K}_{n-1}^{[1]}} h_1 \left( \frac{k_1}{n} \right) h_2 \left( \frac{k_2}{n} \right) \int_{\partial_{k_1}^- \underline{K}_n} \frac{\partial \tilde{\varphi}_n}{\partial x_{k_2}} (\underline{x}) (\Xi_{n,a,b}^n)^{-1} q_n^{0,n}(\underline{x}) \\
        & \qquad \qquad \qquad \qquad \qquad \qquad \qquad \qquad \mathrm{d} \underline{x}_{[0,k_1/n]}^{a,f^-(k_1/n)} \ \mathrm{d} \underline{x}_{[k_1/n,1]}^{f^-(k_1/n), b} \\
        &= \sum_{(k_1, k_2) \in \mathcal{K}_{n-1}^{[2]}} h_1 \left( \frac{k_1}{n} \right) h_2 \left( \frac{k_2}{n} \right) \int_{\partial_{k_1}^- \underline{K}_n} \frac{\partial \tilde{\varphi}_n}{\partial x_{k_2}} (\underline{x}) (\Xi_{n,a,b}^n)^{-1} q_n^{0,n}(\underline{x}) \\
        & \qquad \qquad \qquad \qquad \qquad \qquad \qquad \qquad \mathrm{d} \underline{x}_{[0,k_1/n]}^{a,f^-(k_1/n)} \ \mathrm{d} \underline{x}_{[k_1/n,1]}^{f^-(k_1/n), b} \\
        &\quad + \sum_{(k_2, k_1) \in \mathcal{K}_{n-1}^{[2]}} h_1 \left( \frac{k_1}{n} \right) h_2 \left( \frac{k_2}{n} \right) \int_{\partial_{k_1}^- \underline{K}_n} \frac{\partial \tilde{\varphi}_n}{\partial x_{k_2}} (\underline{x}) (\Xi_{n,a,b}^n)^{-1} q_n^{0,n}(\underline{x}) \\
        & \qquad \qquad \qquad \qquad \qquad \qquad \qquad \qquad \mathrm{d} \underline{x}_{[0,k_1/n]}^{a,f^-(k_1/n)} \ \mathrm{d} \underline{x}_{[k_1/n,1]}^{f^-(k_1/n), b} \\
        &= - \sum_{(k_1, k_2) \in \mathcal{K}_{n-1}^{[2]}} h_1 \left( \frac{k_1}{n} \right) h_2 \left( \frac{k_2}{n} \right) \int_{\partial_{k_1}^- \underline{K}_n} \tilde{\varphi}_n (\underline{x}) (\Xi_{n,a,b}^n)^{-1} \frac{\partial q_n^{0,n}}{\partial x_{k_2}}(\underline{x}) \\
        & \qquad \qquad \qquad \qquad \qquad \qquad \qquad \qquad \mathrm{d} \underline{x}_{[0,k_1/n]}^{a,f^-(k_1/n)} \ \mathrm{d} \underline{x}_{[k_1/n,1]}^{f^-(k_1/n), b} \\
        &\quad + \sum_{(k_1, k_2) \in \mathcal{K}_{n-1}^{[2]}} h_1 \left( \frac{k_1}{n} \right) h_2 \left( \frac{k_2}{n} \right) \int_{\partial_{k_2}^+ \partial_{k_1}^- \underline{K}_n} \tilde{\varphi}_n (\underline{x}) (\Xi_{n,a,b}^n)^{-1} q_n^{0,n}(\underline{x}) \\
        &\qquad \qquad \qquad \qquad \times \ \mathrm{d} \underline{x}_{[0,k_1/n]}^{a,f^-(k_1/n)} \ \mathrm{d} \underline{x}_{[k_1/n, k_2/n]}^{f^-(k_1/n), f^+(k_2/n)} \ \mathrm{d} \underline{x}_{[k_2/n,1]}^{f^+(k_2/n), b} \\
        &\quad - \sum_{(k_1, k_2) \in \mathcal{K}_{n-1}^{[2]}} h_1 \left( \frac{k_1}{n} \right) h_2 \left( \frac{k_2}{n} \right) \int_{\partial_{k_2}^- \partial_{k_1}^- \underline{K}_n} \tilde{\varphi}_n (\underline{x}) (\Xi_{n,a,b}^n)^{-1} q_n^{0,n}(\underline{x}) \\
        &\qquad \qquad \qquad \qquad \times \ \mathrm{d} \underline{x}_{[0,k_1/n]}^{a,f^-(k_1/n)} \ \mathrm{d} \underline{x}_{[k_1/n, k_2/n]}^{f^-(k_1/n), f^-(k_2/n)} \ \mathrm{d} \underline{x}_{[k_1/n,1]}^{f^-(k_2/n), b} \\
        &\quad - \sum_{(k_2, k_1) \in \mathcal{K}_{n-1}^{[2]}} h_1 \left( \frac{k_1}{n} \right) h_2 \left( \frac{k_2}{n} \right) \int_{\partial_{k_1}^- \underline{K}_n} \tilde{\varphi}_n (\underline{x}) (\Xi_{n,a,b}^n)^{-1} \frac{\partial q_n^{0,n}}{\partial x_{k_2}}(\underline{x}) \\
        &\qquad \qquad \qquad \qquad \qquad \qquad \qquad \qquad \mathrm{d} \underline{x}_{[0,k_1/n]}^{a,f^-(k_1/n)} \ \mathrm{d} \underline{x}_{[k_1/n,1]}^{f^-(k_1/n), b} \\
        &\quad + \sum_{(k_2, k_1) \in \mathcal{K}_{n-1}^{[2]}} h_1 \left( \frac{k_1}{n} \right) h_2 \left( \frac{k_2}{n} \right) \int_{\partial_{k_2}^+ \partial_{k_1}^- \underline{K}_n} \tilde{\varphi}_n (\underline{x}) (\Xi_{n,a,b}^n)^{-1} q_n^{0,n}(\underline{x}) \\
        &\qquad \qquad \qquad \qquad \times \ \mathrm{d} \underline{x}_{[0,k_2/n]}^{a,f^+(k_2/n)} \ \mathrm{d} \underline{x}_{[k_2/n, k_1/n]}^{f^+(k_2/n), f^-(k_1/n)} \ \mathrm{d} \underline{x}_{[k_1/n,1]}^{f^-(k_1/n), b} \\
        &\quad - \sum_{(k_2, k_1) \in \mathcal{K}_{n-1}^{[2]}} h_1 \left( \frac{k_1}{n} \right) h_2 \left( \frac{k_2}{n} \right) \int_{\partial_{k_2}^- \partial_{k_1}^- \underline{K}_n} \tilde{\varphi}_n (\underline{x}) (\Xi_{n,a,b}^n)^{-1} q_n^{0,n}(\underline{x}) \\
        &\qquad \qquad \qquad \qquad \times \ \mathrm{d} \underline{x}_{[0,k_2/n]}^{a,f^-(k_2/n)} \ \mathrm{d} \underline{x}_{[k_2/n, k_1/n]}^{f^-(k_2/n), f^-(k_1/n)} \ \mathrm{d} \underline{x}_{[k_1/n,1]}^{f^-(k_1/n), b}
    \end{aligned}
\end{equation}

\begin{equation}\nonumber
    \begin{aligned}
        &= - \sum_{k_1, k_2 \in \mathcal{K}_{n-1}^{[1]}} h_1 \left( \frac{k_1}{n} \right) h_2 \left( \frac{k_2}{n} \right) \int_{\partial_{k_1}^- \underline{K}_n} \tilde{\varphi}_n (\underline{x}) (\Xi_{n,a,b}^n)^{-1} \frac{\partial q_n^{0,n}}{\partial x_{k_2}}(\underline{x}) \\
        &\qquad \qquad \qquad \qquad \qquad \qquad \qquad \qquad \mathrm{d} \underline{x}_{[0,k_1/n]}^{a,f^-(k_1/n)} \ \mathrm{d} \underline{x}_{[k_1/n,1]}^{f^-(k_1/n), b} \\
        &\quad + \sum_{(k_1, k_2) \in \mathcal{K}_{n-1}^{[2]}} h_1 \left( \frac{k_1}{n} \right) h_2 \left( \frac{k_2}{n} \right) \int_{\partial_{k_2}^+ \partial_{k_1}^- \underline{K}_n} \tilde{\varphi}_n (\underline{x}) (\Xi_{n,a,b}^n)^{-1} q_n^{0,n}(\underline{x}) \\
        &\qquad \qquad \qquad \qquad \times \ \mathrm{d} \underline{x}_{[0,k_1/n]}^{a,f^-(k_1/n)} \ \mathrm{d} \underline{x}_{[k_1/n, k_2/n]}^{f^-(k_1/n), f^+(k_2/n)} \ \mathrm{d} \underline{x}_{[k_2/n,1]}^{f^+(k_2/n), b} \\
        &\quad + \sum_{(k_2, k_1) \in \mathcal{K}_{n-1}^{[2]}} h_1 \left( \frac{k_1}{n} \right) h_2 \left( \frac{k_2}{n} \right) \int_{\partial_{k_2}^+ \partial_{k_1}^- \underline{K}_n} \tilde{\varphi}_n (\underline{x}) (\Xi_{n,a,b}^n)^{-1} q_n^{0,n}(\underline{x}) \\
        &\qquad \qquad \qquad \qquad \times \ \mathrm{d} \underline{x}_{[0,k_2/n]}^{a,f^+(k_2/n)} \ \mathrm{d} \underline{x}_{[k_2/n, k_1/n]}^{f^+(k_2/n), f^-(k_1/n)} \ \mathrm{d} \underline{x}_{[k_1/n,1]}^{f^-(k_1/n), b} \\
        &\quad - \sum_{(k_1, k_2) \in \mathcal{K}_{n-1}^{[2]}} \left( h_1 \left( \frac{k_1}{n} \right) h_2 \left( \frac{k_2}{n} \right) + h_1 \left( \frac{k_2}{n} \right) h_2 \left( \frac{k_1}{n} \right) \right) \\
        &\qquad \qquad \qquad \qquad \int_{\partial_{k_2}^- \partial_{k_1}^- \underline{K}_n} \tilde{\varphi}_n (\underline{x}) (\Xi_{n,a,b}^n)^{-1} q_n^{0,n}(\underline{x}) \\
        &\qquad \qquad \qquad \qquad \times \ \mathrm{d} \underline{x}_{[0,k_1/n]}^{a,f^-(k_1/n)} \ \mathrm{d} \underline{x}_{[k_1/n, k_2/n]}^{f^-(k_1/n), f^-(k_2/n)} \ \mathrm{d} \underline{x}_{[k_2/n,1]}^{f^-(k_2/n), b}
    \end{aligned}
\end{equation}

Hence, we obtain
\begin{equation}\nonumber
    \begin{aligned}
        &-I_n + B_n^+ - B_n^- \\
        &= \sum_{k_1, k_2 \in \mathcal{K}_{n-1}^{[1]}} h_1 \left( \frac{k_1}{n} \right) h_2 \left( \frac{k_2}{n} \right) \int_{\underline{K}_n} \tilde{\varphi}_n (\underline{x}) (\Xi_{n,a,b}^n)^{-1} \frac{\partial^2 q_n^{0,n}}{\partial x_{k_2} \partial x_{k_1}}(\underline{x}) \ \mathrm{d} \underline{x}^{a,b} \\
        &\quad - \sum_{k_1, k_2 \in \mathcal{K}_{n-1}^{[1]}} h_1 \left( \frac{k_1}{n} \right) h_2 \left( \frac{k_2}{n} \right) \int_{\partial_{k_1}^+ \underline{K}_n} \tilde{\varphi}_n (\underline{x}) (\Xi_{n,a,b}^n)^{-1} \frac{\partial q_n^{0,n}}{\partial x_{k_2}}(\underline{x}) \\
        &\qquad \qquad \qquad \qquad \qquad \qquad \qquad \qquad \mathrm{d} \underline{x}_{[0,k_1/n]}^{a,f^+(k_1/n)} \ \mathrm{d} \underline{x}_{[k_1/n,1]}^{f^+(k_1/n), b} \\
        &\quad + \sum_{k_1, k_2 \in \mathcal{K}_{n-1}^{[1]}} h_1 \left( \frac{k_1}{n} \right) h_2 \left( \frac{k_2}{n} \right) \int_{\partial_{k_1}^- \underline{K}_n} \tilde{\varphi}_n (\underline{x}) (\Xi_{n,a,b}^n)^{-1} \frac{\partial q_n^{0,n}}{\partial x_{k_2}}(\underline{x}) \\
        &\qquad \qquad \qquad \qquad \qquad \qquad \qquad \qquad \mathrm{d} \underline{x}_{[0,k_1/n]}^{a,f^-(k_1/n)} \ \mathrm{d} \underline{x}_{[k_1/n,1]}^{f^-(k_1/n), b} \\
        &\quad + \sum_{(k_1, k_2) \in \mathcal{K}_{n-1}^{[2]}} \left( h_1 \left( \frac{k_1}{n} \right) h_2 \left( \frac{k_2}{n} \right) + h_1 \left( \frac{k_2}{n} \right) h_2 \left( \frac{k_1}{n} \right) \right) \\
        &\qquad \qquad \qquad \qquad \times \int_{\partial_{k_2}^+ \partial_{k_1}^+ \underline{K}_n} \tilde{\varphi}_n (\underline{x}) (\Xi_{n,a,b}^n)^{-1} q_n^{0,n}(\underline{x}) \\
        &\qquad \qquad \qquad \qquad \times \ \mathrm{d} \underline{x}_{[0,k_1/n]}^{a,f^+(k_1/n)} \ \mathrm{d} \underline{x}_{[k_1/n, k_2/n]}^{f^+(k_1/n), f^+(k_2/n)} \ \mathrm{d} \underline{x}_{[k_2/n,1]}^{f^+(k_2/n), b} \\
        &\quad - \sum_{(k_1, k_2) \in \mathcal{K}_{n-1}^{[2]}} \left( h_1 \left( \frac{k_1}{n} \right) h_2 \left( \frac{k_2}{n} \right) + h_1 \left( \frac{k_2}{n} \right) h_2 \left( \frac{k_1}{n} \right) \right) \\
        &\qquad \qquad \qquad \qquad \times \int_{\partial_{k_2}^- \partial_{k_1}^+ \underline{K}_n} \tilde{\varphi}_n (\underline{x}) (\Xi_{n,a,b}^n)^{-1} q_n^{0,n}(\underline{x}) \\
        &\qquad \qquad \qquad \qquad \times \ \mathrm{d} \underline{x}_{[0,k_1/n]}^{a,f^+(k_1/n)} \ \mathrm{d} \underline{x}_{[k_1/n, k_2/n]}^{f^+(k_1/n), f^-(k_2/n)} \ \mathrm{d} \underline{x}_{[k_1/n,1]}^{f^-(k_2/n), b} \\
        &\quad - \sum_{(k_1, k_2) \in \mathcal{K}_{n-1}^{[2]}} \left( h_1 \left( \frac{k_1}{n} \right) h_2 \left( \frac{k_2}{n} \right) + h_1 \left( \frac{k_2}{n} \right) h_2 \left( \frac{k_1}{n} \right) \right) \\
        &\qquad \qquad \qquad \qquad \times \int_{\partial_{k_2}^+ \partial_{k_1}^- \underline{K}_n} \tilde{\varphi}_n (\underline{x}) (\Xi_{n,a,b}^n)^{-1} q_n^{0,n}(\underline{x}) \\
        &\qquad \qquad \qquad \qquad \times \ \mathrm{d} \underline{x}_{[0,k_1/n]}^{a,f^-(k_1/n)} \ \mathrm{d} \underline{x}_{[k_1/n, k_2/n]}^{f^-(k_1/n), f^+(k_2/n)} \ \mathrm{d} \underline{x}_{[k_2/n,1]}^{f^+(k_2/n), b} \\
        &\quad + \sum_{(k_1, k_2) \in \mathcal{K}_{n-1}^{[2]}} \left( h_1 \left( \frac{k_1}{n} \right) h_2 \left( \frac{k_2}{n} \right) + h_1 \left( \frac{k_2}{n} \right) h_2 \left( \frac{k_1}{n} \right) \right) \\
        &\qquad \qquad \qquad \qquad \times \int_{\partial_{k_2}^- \partial_{k_1}^- \underline{K}_n} \tilde{\varphi}_n (\underline{x}) (\Xi_{n,a,b}^n)^{-1} q_n^{0,n}(\underline{x}) \\
        &\qquad \qquad \qquad \qquad \times \ \mathrm{d} \underline{x}_{[0,k_1/n]}^{a,f^-(k_1/n)} \ \mathrm{d} \underline{x}_{[k_1/n, k_2/n]}^{f^-(k_1/n), f^-(k_2/n)} \ \mathrm{d} \underline{x}_{[k_2/n,1]}^{f^-(k_2/n), b} .
    \end{aligned}
\end{equation}
Therefore, the integral by parts formula is obtained by taking the limit.

\end{appendices}

\newpage
\begin{flushleft}
\mbox{  }\\
\hspace{95mm} Kensuke Ishitani\\
\hspace{95mm} Department of Mathematical Sciences\\
\hspace{95mm} Tokyo Metropolitan University\\
\hspace{95mm} Hachioji, Tokyo 192-0397\\
\hspace{95mm} Japan\\
\hspace{95mm} e-mail: k-ishitani@tmu.ac.jp\\
\mbox{  }\\
\hspace{95mm} Soma Nishino\\
\hspace{95mm} Department of Mathematical Sciences\\
\hspace{95mm} Tokyo Metropolitan University\\
\hspace{95mm} Hachioji, Tokyo 192-0397\\
\hspace{95mm} Japan\\
\hspace{95mm} e-mail: nishino-soma@ed.tmu.ac.jp\\
\end{flushleft}

\end{document}